\documentclass[smallextended]{svjour3}

\usepackage{graphicx}
\usepackage{psfrag} 
\usepackage{mathptmx}      
\usepackage{amssymb}
\usepackage{amsmath}
\usepackage{algorithm}
\usepackage{marvosym}
\usepackage{a4wide}
\usepackage[noend]{algpseudocode}
\usepackage{xspace}
\algrenewcommand\algorithmicrequire{\textbf{Input:}}
\algrenewcommand\algorithmicensure{\textbf{Output:}}

\usepackage{color}

\usepackage[mathscr]{euscript}

\newtheorem{observation}[theorem]{Observation}

\newcommand{\lca}{\ensuremath{\operatorname{lca}}}

\newcommand{\gatex}{\textsc{GaTEx}\xspace}

\newcommand{\new}[1]{\begingroup\color{black}#1\endgroup}

\begin{document}

\title{Representing distance-hereditary graphs with multi-rooted trees}

\author{Guillaume E. Scholz}

\institute{G. E. Scholz \at Bioinformatics Group, Department of Computer Science and Interdisciplinary Center for Bioinformatics, Universit\"at
				Leipzig, Leipzig, Germany.
              \email{guillaume@bioinf.uni-leipzig.de}}

\date{Received: date / Accepted: date}

\maketitle

\begin{abstract}
Arboreal networks are a generalization of rooted trees, defined by keeping the tree-like structure, but dropping the requirement for a single root. Just as the class of cographs is precisely the class of undirected graphs that can be explained by a labelled rooted tree $(T,t)$, we show that the class of distance-hereditary graphs is precisely the class of undirected graphs that can be explained by a labelled arboreal network $(N,t)$.
\end{abstract}

\section{Introduction}

Among other interesting properties, cographs \cite{Corneil:81,Seinsche74,Sumner74} are precisely the undirected graphs $G$ that can be represented by a labelled tree, that is, a pair $(T,t)$, where $T$ is a rooted tree with leaf set $V(G)$, and $t$ is a labelling of the vertices of $T$ into the set $\{0,1\}$. Specifically, we say that the pair $(T,t)$ explains $G$ if, for any two distinct vertices $x$ and $y$ of $G$, $x$ and $y$ are joined by an edge in $G$ if and only if the least common ancestor of $x$ and $y$ in $T$ has label $1$ via $t$.

Such a property is of particular interest in computer science, and many problems that are known to be NP-hard for general classes of graphs can be solved in polynomial time for cographs \cite{BLS99,Corneil:85,Gao:13}, using the labelled tree $(T,t)$ as a guide. The relation between cographs and labelled rooted trees is also exploited in evolutionary biology \cite{HHH13,HW:16b,H:15a,HHH12,LM16,LM14,LM15}, a field in which labelled rooted trees play a key role in understanding the common evolutionary history of a set of taxa. For example, given a set of genes $X$, it is possible to infer the pairs of genes $\{x,y\}$ of $X$ which descend from a speciation event \cite{AD09}. This information can be stored in an undirected graph $G$ with vertex set $X$. If that graph $G$ is a cograph, a labelled tree $(T,t)$ explaining $G$ is therefore a good candidate for describing the common evolutionary history of genes in $X$.

Thus, the interest of extending the notion of graph representations beyond the pair cograph-rooted tree is (at least) twofold. First, identify further graph classes for which NP-hard problems become polynomial solvable, and design algorithm to solve such problems. Second, provide evolutionary biologists with tools to deal with data that do not fit within the existing framework. In the context of orthology relation, this means providing a way to represent the common evolutionary history of a set of genes $X$, when the graph $G$ inferred from this set of genes is not a cograph.

This quest for generalization can be summarized as follows. Given an undirected graph $G$ with vertex set $X$, is there a pair $(N,t)$ such that $N$ is a directed graph with leaf set $X$, and $t$ is a labelling of the vertices of $N$ into the set $\{0,1\}$, such that $x$ and $y$ are joined by an edge in $G$ if and only if the least common ancestor of $x$ and $y$ in $T$ has label $1$ via $t$? Naturally, such a definition requires the notion of a least common ancestor in $N$ to be defined without ambiguity.

Preliminary work into this question \cite{BHS21} showed that there exists a class of directed graphs, known as rooted median graphs, for which the notion of least common ancestor is well-defined, and is such that any graph $G$ can be explained by a pair $(N,t)$ such that $N$ is a rooted median graph. The rooted median graphs constructed in \cite{BHS21}, however, have $O(|V(G)|^2)$ vertices, which make them both impractical for computational purposes, and unlikely to represent a plausible evolutionary scenario.

In \cite{HS22} and \cite{HS23}, the question was restricted to a special class of rooted networks $N$, known as galled-trees \cite{G03,HSS22,RV09}, a class often considered to be "one step away" from rooted trees in terms of complexity \cite{GBP12,GHK17,HS18}. This led to the definition of a new class of undirected graphs, \gatex graphs, and to the design of polynomial-time algorithms to solve some of the classical optimization problems on graphs from that particular class \cite{HS23b}.

In this contribution, we turn our attention to a newly introduced class of directed graphs, the class of arboreal networks \cite{HMS24,HMS22}. Such structures are treelike, in the sense that their underlying, unidrected graph is a tree. However, they differ from rooted tree in the fact that they are allowed to have more than one single root. In evolutionary studies, dropping the requirement for a single root may be useful, for example when the ancestral relationship between distinct lineages is too ancient to be inferred with accuracy \cite{scholz19}. In the main result of this contribution, Theorem~\ref{th:dh}, we show that the class of undirected graphs that can be explained be a labelled arboreal network is precisely the class of distance-hereditary graphs, a well-known and well studied class of undirected graphs \cite{DiS:12,GP12,HM:90,H77}.

This contribution is organized as follows. In Section~\ref{sec-ug}, we introduce the basic definitions and notations regarding undirected graphs. In addition, we define three key classes of undirected graphs, cographs, Ptolemaic graphs and distance-hereditary graphs, and we present some important properties of these three classes of graphs. In Section~\ref{sec-ntw}, we define the concepts of networks, arboreal networks and labelled networks, and we formally define the notion of an arboreal-explainable graph, that is, an undirected graph that can be explained by a labelled arboreal network. In Section~\ref{sec-con}, we present some  useful properties of arboreal-explainable graphs regarding connectedness. More specifically, we show in Proposition~\ref{pr:red} that a graph $G$ is arboreal-explainable if and only if all connected components of $G$ are arboreal-explainable.

In Section~\ref{sec-char}, we present our main result, Theorem~\ref{th:dh}, which states that arboreal-explainable graphs are precisely the distance-hereditary graphs. In addition, we provide an algorithm, Algorithm~\ref{alg:dh}, which constructs, for a given distance-hereditary graph $G$, a labelled arboreal network $(N,t)$ explaining $G$.

Interestingly, arboreal networks, closely related to the laminar trees introduced in \cite{UU09}, have been shown to bear strong links to Ptolemaic graphs \cite{HMS24,UU09}, a subclass of distance-hereditary graphs \cite{HM:90,H81}. More specifically, a given arboreal network $N$ is associated to a unique Ptolemaic graph $G^*$ with vertex set $L(N)$, which records the pairs of leaves of $N$ that have a common ancestor in $N$. In section~\ref{sec-compl}, we investigate the relation between the graph $G$ explained by a given labelled arboreal network $(N,t)$, and the Ptolemaic graph $G^*$ associated to $N$. More specifically, we characterize the pairs $(G,G^*)$ of undirected graphs for which there exists a labelled arboreal network $(N,t)$ explaining $G$, such that $G^*$ is the Ptolemaic graph associated to $N$ (Theorem~\ref{th:ptext}). Finally, in Section ~\ref{sec-2r}, we restrict our attention to labelled arboreal networks with exactly two roots, and we charactarize the subclass of arboreal-explainable graphs that can be explained by such a network (Theorem~\ref{th:2rgt}). We conclude in Section~\ref{sec-out} by listing some open question raised by the findings presented in this contribution.

\section{Three classes of undirected graphs}\label{sec-ug}

An \emph{undirected graph} $G$ is a pair $(X,E)$, where $X$ is a nonempty finite set, and $E$ is a set of unordered pairs of distinct elements of $X$. We call $X$ the \emph{vertex set} of $G$, $E$ the \emph{edge set} of $G$. We also sometimes denote the vertex set of $G$ by $V(G)$, and the edge set of $G$ by $E(G)$. If $x,y \in X$ are such that $\{x,y\} \in E$, we say that $x$ and $y$ are \emph{joined by an edge} in $G$.

Let $G=(X,E)$ be an undirected graph. We say that a graph $H=(X',E')$ is a \emph{subgraph} of $G$ if $X' \subseteq X$ and $E' \subseteq E$, and we say that $H$ is an \emph{induced subgraph} of $G$ if $X' \subseteq X$ and $E'=\{\{u,v\} \in E, u \in X', v \in X'\}$. In that case, we write $H=G[X']$.

An undirected graph $G=(X,E)$ is called a \emph{clique} if $E$ is the set of all possible pairs of distinct elements of $X$. If there exists an ordering $x_1, \ldots, x_k$, $k=|X|$ of the elements of $X$ such that $E=\{\{x_i,x_{i+1}\}, i \in \{1, \ldots, k-1\}\}$, we call $G$ a \emph{path}. In that case, we say that $G$ is a path between $x_1$ and $x_k$, and \new{we say that the \emph{length} of $G$ is $k-1$}. Similarly, we say that $G$ is a \emph{cycle} if there exists an ordering $x_1, \ldots, x_k$, $k=|X|$ of the elements of $X$ such that \new{$E=\{\{x_i,x_{i+1}\}, i \in \{1, \ldots, k-1\}\} \cup \{\{x_k,x_1\}\}$}. In that case, we call $k$ the \emph{size} of $G$.

We say that two vertices $x,y$ of an undirected graph $G$ are \emph{connected} if there exists a path in $G$ between $x$ and $y$. In that case, we define $d_G(x,y)$ as the length of the shortest path between $x$ and $y$ in $G$. If any two vertices $x,y$ of $G$ are connected, we say that $G$ is \emph{connected}. If \new{$|X| \geq 3$} and for all $x \in V(G)$, the graph $G[X \setminus \{x\}]$ is connected, we say that $G$ is \emph{biconnected}. For $Y \subseteq X$, we say that $G[Y]$ is a \emph{connected component} (\emph{resp.} a \emph{biconnected component}) of $G$ if $G[Y]$ is connected (\emph{resp.} biconnected), and no proper superset $Y'$ of $Y$ is such that $G[Y']$ is connected (\emph{resp.} biconnected). Finally, we say that $G$ is a \emph{tree} if it does not contain an induced cycle, \emph{chordal} if it does not contain an induced cycle of size $4$ or more, and \emph{hole-free} if it does not contain an induced cycle of size $5$ or more.

We next define two key operations on undirected graphs. Let $G_1=(X_1,E_1)$ and $G_2=(X_2,E_2)$ be two nonempty graphs with $X_1 \cap X_2=\emptyset$. The \emph{disjoint union} of $G_1$ and $G_2$ is the graph with vertex set $X_1 \cup X_2$ and edge set $E_1 \cup E_2$. The \emph{joint} of $G_1$ and $G_2$ is the graph with vertex set $X_1 \cup X_2$ and edge set $E_1 \cup E_2 \cup \{\{u,v\}, u \in X_1, v \in X_2\}$.

The disjoint union and joint operations are at the heart of the definition of \emph{cographs} \cite{Corneil:81,Seinsche74,Sumner74}. Cographs are defined recursively, using the following three rules:
\begin{itemize}
\item[(C1)] A single-vertex graph is a cograph.
\item[(C2)] The disjoint union of two cographs is a cograph.
\item[(C3)] The joint of two cographs is a cograph.
\end{itemize}

Interestingly, cographs are precisely the graphs that do not contain an induced path of length $3$ (or a $P_4$, for short) as an induced subgraph. This characterization will prove useful in the remaining of this paper.

An undirected graph $G=(X,E)$ is \emph{Ptolemaic} \cite{H81,KC65} if for all $x,y,z,t \in X$ distinct and pairwise connected, $d_G(x,y)d_G(z,t)+d_G(x,t)d_G(y,z) \geq d_G(x,z)d_G(y,t)$. The latter inequality is known as \emph{Ptolemy's inequality}.

\begin{figure}[h]
\begin{center}
\includegraphics[scale=0.8]{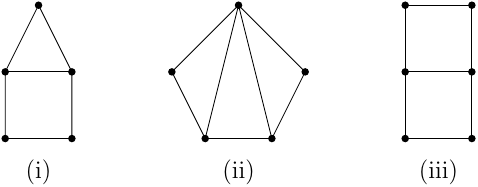}
\end{center}
\caption{(i) The house. (ii) The gem. (iii) The domino.}
\label{fig-forb}
\end{figure}

An undirected graph $G$ is \emph{distance-hereditary} \cite{H77}, if for all pair of connected vertices $x,y$ of $G$, and for all connected induced subgraph $H$ of $G$ containing $x$ and $y$, $d_G(x,y)=d_H(x,y)$ holds. It has been shown in \cite{HM:90} (see also \cite{DiS:12}) that a graph $G$ is distance-hereditary if and only if $G$ is hole-free, and does not contain the \emph{house}, the \emph{gem} and the \emph{domino} as an induced subgraph (see Figure~\ref{fig-forb} for a depiction of these three graphs). In particular, all graphs with $4$ or less vertices are distance-hereditary.

Equivalently (\cite[Theorem 1]{BM86}), a graph $G$ is distance-hereditary if it can be recursively constructed from a single vertex by the following operations:
\begin{itemize}
\item[(R1)] Add a \emph{pendant-vertex}: pick a vertex $x$. Add a new vertex $y$ and the edge $\{x,y\}$.
\item[(R2)] Create \emph{false-twins}: pick a vertex $x$. Add a new vertex $y$ and add the edge $\{y,z\}$ for all vertices $z$ adjacent to $x$.
\item[(R3)] Create \emph{true-twins}: pick a vertex $x$. Add a new vertex $y$, add the edge $\{x,y\}$, and add the edge $\{y,z\}$ for all vertices $z$ adjacent to $x$.
\end{itemize}

As a direct consequence of this definition, we have:

\begin{lemma}\label{lm:indh}
Let $G=(X,E)$ be an undirected graph. If there exists $x,y \in X$ distinct such that $G[X \setminus \{y\}]$ is distance-hereditary, and $G$ can be obtained from $G[X \setminus \{y\}]$ by adding $y$ as a pendant-vertex, false-twin or true-twin to $x$, then $G$ is distance-hereditary.
\end{lemma}

Interestingly, if only operations (R2) and (R3) are allowed, then the resulting graph $G$ is a cograph. If all three operations are allowed, but operation (R2) is amended with the extra-requirement that the set of vertices $z$ adjacent to $x$ induce a clique, then the resulting graph $G$ is Ptolemaic \cite[Corollary 6]{BM86}. In particular, cographs and Ptolemaic graphs are distance-hereditary.

By \cite{H81}, a graph $G$ is Ptolemaic if and only if $G$ is chordal and distance-hereditary. Note that if $G$ is a chordal graph then $G$ is hole-free, and does not contain the house or the domino as an induced subgraph. Hence, a graph is Ptolemaic if and only if it is chordal, and does not contain the gem as an induced subgraph.

Finally, we say that a graph property $(\mathcal P)$ is \emph{hereditary} if, for all graphs $G$ satisfying property $(\mathcal P)$, it holds that all induced subgraphs $H$ of $G$ satisfy property $(\mathcal P)$. As we have just seen, cographs, Ptolemaic graphs and distance-hereditary graphs can be characterized in terms of a list of forbidden induced-subgraphs. As a consequence, the properties of being a cograph, being Ptolemaic, and being distance-hereditary, are all hereditary.

\section{Labelled arboreal network}\label{sec-ntw}

In this section, we turns our attention to directed graphs. A directed graph $G$ is a pair $(X,E)$ where $X$ is a nonempty finite set, and $E$ is a set of ordered pairs of distinct elements of $X$. When applicable, all definitions and terminology are similar to those introduced in the previous section, with the exception of the name given to the elements of $E$, which we call \emph{arcs} in the directed case. For $(u,v) \in E$, we say that $u$ is the \emph{parent} of $v$ (in $G$), and $v$ is the \emph{child} of $u$ (in $G$). \new{Given a vertex $v$, we call the \emph{indegree} of $v$ the number of parents of $v$, the \emph{outdegree} of $v$ the number of children of $v$, and the \emph{degree} of $v$ the sum of its indegree and its outdegree.} For $G=(X,E)$ a directed graph, the \emph{underlying undirected graph} of $G$ is the graph $\overline G=(X,\overline E)$, where $\overline E$ is the set of all pairs $\{u,v\}$ of elements of $X$ such that (at least) one of $(u,v)$ or $(v,u)$ is an element of $E$. Informally speaking, $\overline G$ is the undirected graph obtained from $G$ by ignoring the direction of the arcs. \new{If $G$ does not contain a pair $u,v$ of vertices such that both $(u,v)$ and $(v,u)$ are arcs of $G$, then the degree of a vertex $v$ in $G$ coincides with the degree of $v$ in $\overline G$.}

\new{We say that a directed graph $G=(X,E)$ is a \emph{(directed) path} if there exists an ordering $x_1, \ldots, x_k$, $k=|X|$ of the elements of $X$ such that $E=\{(x_i,x_{i+1}), i \in \{1, \ldots, k-1\}\}$. In that case, we say that $G$ is a path from $x_1$ to $x_k$, and we say that the \emph{length} of $G$ is $k-1$ (that is, the number of edges of $G$). Similarly, we say that $G$ is a \emph{directed cycle} if there exists an ordering $x_1, \ldots, x_k$, $k=|X|$ of the elements of $X$ such that $E=\{(x_i,x_{i+1}), i \in \{1, \ldots, k-1\}\} \cup \{(x_k,x_1)\}$. A directed graph that does not contain any directed cycle is called \emph{acyclic}. Finally, we say that $G$ is \emph{connected} if the underlying undirected graph $\overline G$ of $G$ is connected.}

Two key operations on directed graphs will be useful throughout this paper. For $G=(X,E)$ a directed graph, and $(u,v) \in E$, \emph{subdividing} $(u,v)$ is the operation of removing $(u,v)$ from $E$, adding a new vertex $w$ to $X$, and adding the arcs $(u,w)$ and $(w,v)$ to $E$. Note that, in the resulting graph, $w$ has indegree and outdegree $1$. The reverse operation is the \emph{suppression} of a vertex $w$ of indegree $1$ and outdegree $1$. It consists in removing $w$ and its incident arcs, and add a new arc from the parent of $w$ to the child of $w$.

\begin{figure}[h]
\begin{center}
\includegraphics[scale=0.8]{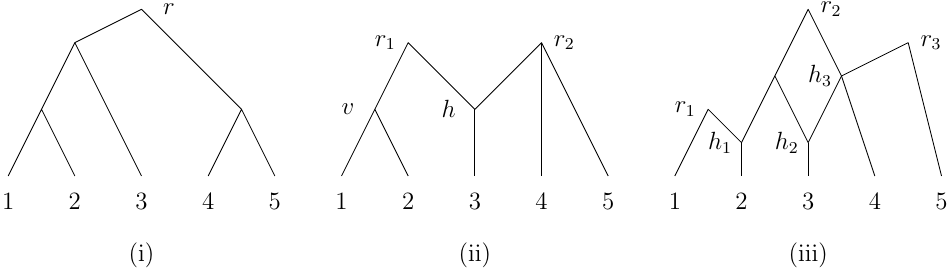}
\end{center}
\caption{(i) A phylogenetic tree with leaf set $X=\{1,2,3,4,5\}$. (ii) An arboreal network with two roots and leaf set $X$. A network with three roots and leaf set $X$ that is not arboreal. Here and in all subsequent Figures depicting networks, all arcs are assumed to be directed downwards. Note that only the network in (i) is binary.}
\label{fig-ntw}
\end{figure}

Following \cite{HMS24}, we say that a connected, directed acyclic graph $N$ is a \emph{network} if all vertices of $N$ of indegree $0$ have outdegree $2$ or more, all vertices of $N$ of outdegree $0$ have indegree $1$, and no vertex of $N$ has indegree and outdegree $1$.  We say that $N$ is \emph{arboreal} if the underlying undirected graph $\overline N$ of $N$ is a tree. Equivalently, $N$ is arboreal if for all arcs $a$ of $N$, the removal of $a$ from $N$ disconnects $N$. We say that $N$ is a \emph{phylogenetic network} if $N$ has exactly one vertex of indegree $0$, and a \emph{phylogenetic tree} if $N$ is an arboreal phylogenetic network. Figure~\ref{fig-ntw} depicts a phylogenetic tree (i), an arboreal network (ii) and a network that is not a phylogenetic network nor an arboreal network (iii).

We say that a vertex $v$ of $N$ is a \emph{root} if $v$ has indegree $0$, and a \emph{leaf} if $v$ has outdegree $0$. We denote by $R(N)$ the set of roots of $N$, by $H(N)$ the set of vertices of $N$ of indegree $2$ or more, by $V^*(N)$ the set of vertices of $N$ of outdegree $2$ or more, and by $L(N)$ the set of leaves of $N$. For example, for $N$ the network depicted in Figure~\ref{fig-ntw}~(ii), we have $L(N)=\{1,2,3,4,5\}$, $R(N)=\{r_1,r_2\}$, $H(N)=\{h\}$, and $V^*(N)=\{r_1,r_2,v\}$. We say that $N$ is \emph{binary} if all vertices of $V^*(N)$ have outdegree $2$ and indegree $0$ or $1$, and all vertices of $H(N)$ have indegree $2$ and outdegree $1$. In particular, if $N$ is binary, $V^*(N) \cap H(N)=\emptyset$.

Interestingly, the sets $R(N)$ and $H(N)$ are closely related, even more so if $N$ is arboreal:

\begin{lemma}[\cite{HMS24}, Lemma~3.1 and Proposition~3.2]\label{lm:rh}
Let $N$ be a network, $r=|R(N)|$, and $\tilde h=\sum\limits_{h \in H(N)} \mathrm{indeg}(h)-1$. Then, $\tilde h \geq r-1$, and the equality holds if and only if $N$ is arboreal.
\end{lemma}

Note that Lemma~\ref{lm:rh} together with the definition of a phylogenetic tree, immediately implies that a phylogenetic network $N$ is a phylogenetic tree if and only if $H(N)=\emptyset$, which is consistent with the standard definition of a phylogenetic tree.

\new{An important property of arboreal networks is given by the following.

\begin{lemma}\label{lm:uniqp}
Let $N$ be an arboreal network, and let $u,v$ be two vertices of $N$. If $u$ is an ancestor of $v$, then there exists a unique path from $u$ to $v$ in $N$.
\end{lemma}

\begin{proof}
Suppose for contradiction that there exists two distinct paths $P^1$ and $P^2$ from $u$ to $v$ in $N$. Without loss of generality, we may assume that $u$ and $v$ are the only vertices common to $P^1$ and $P^2$. Then, the graph $G$ formed of the union of $P^1$ and $P^2$ induces a cycle in $\overline N$, so $N$ is not arboreal.
\end{proof}

We remark in passing that the converse of Lemma~\ref{lm:uniqp} is not true, in general. In particular, this uniqueness property cannot be used as a characterization for arboreal networks.}

Continuing with the definitions, for $N$ a network and $u,v$ two vertices of $N$, we say that $u$ is an \emph{ancestor} of $v$ (and $v$ a \emph{descendant} of $u$) if there exists a directed path in $N$ from $u$ to $v$. If $u$ and $v$ are distinct, we call $u$ a \emph{proper ancestor} of $v$ (and $v$ a \emph{proper descendant} of $u$). We say that two distinct \new{vertices} $v$ and $w$ of $N$ \emph{share an ancestor} if there exists a vertex $u$ of $N$ that is an ancestor of both $v$ and $w$. \new{For $x$ and $y$ two leaves of $N$ that share an ancestor in $N$}, we say that $u$ is a \emph{least common ancestor} of $x$ and $y$ in $N$ if $u$ is an ancestor of both $x$ and $y$ in $N$, and no child of $u$ enjoys this property. Note that, in general, the least common ancestor of two leaves $x$ and $y$ is not necessarily unique. However, we have:

\begin{lemma}[\cite{HMS24}, Proposition~7.1]\label{lm:lcau}
If $N$ is an arboreal network, and $x$ and $y$ are two distinct leaves of $N$ that share an ancestor in $N$, then $x$ and $y$ have a unique least common ancestor in $N$.
\end{lemma}

In view of Lemma~\ref{lm:lcau}, for $N$ an arboreal network and $x$ and $y$ two distinct leaves of $N$ that share an ancestor in $N$, we denote by $\lca_N(x,y)$ the unique least common ancestor of $x$ and $y$ in $N$. We have:

\begin{lemma}\label{lm:lcas}
Let $N$ be an arboreal network, and let $u \in V(N)$. There exists $x,y \in L(N)$ such that $u=\lca_N(x,y)$ if and only if $u \in V^*(N)$.
\end{lemma}

\begin{proof}
Suppose first that $u=\lca_N(x,y)$ for some $x,y \in L(N)$ distinct. In particular, $u$ is an ancestor of at least two leaves of $N$, so $u \notin L(N)$. Hence, $u$ has outdegree at least $1$ in $N$. If $u$ has outdegree $1$, then $x$ and $y$ are descendant of the unique child $u'$ of $u$, which contradicts the fact that $u=\lca_N(x,y)$. Hence, $u$ has outdegree at least $2$ in $N$, so $u \in V^*(N)$.

Conversely, suppose that $u \in V^*(N)$, and let $v_1$, $v_2$ be two distinct children of $u$ in $N$. Let $x$ and $y$ be descendants of $v_1$ and $v_2$, respectively. Clearly, $x$ and $y$ are descendants of $u$. Moreover,  there is no child $u'$ of $u$ that is an ancestor of both $x$ and $y$, since \new{this would contradict Lemma~\ref{lm:uniqp}}. Hence, we have $u=\lca_N(x,y)$.
\end{proof}

\begin{figure}[h]
\begin{center}
\includegraphics[scale=0.6]{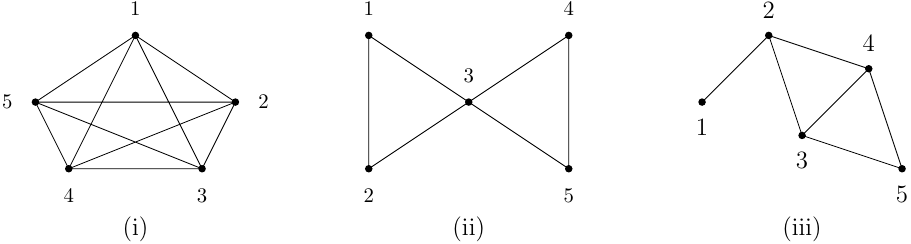}
\end{center}
\caption{(i), (ii), (iii) The shared ancestry graph of the network depicted in Figure~\ref{fig-ntw}~(i), (ii) and (iii), respectively.}
\label{fig-shA}
\end{figure}

For $v$ a vertex of $N$, we denote by $L_v$ the set of all leaves of $N$ that are descendants of $v$. \new{A further key property of arboreal networks is given by the following.

\begin{lemma}\label{lm:emptyint}
Let $N$ be an arboreal network. If $v_1$ and $v_2$ are two incomparable vertices of $N$ that share an ancestor, then $L_{v_1} \cap L_{v_2}=\emptyset$.
\end{lemma}

\begin{proof}
Let $v_1$ and $v_2$ be two incomparable vertices of $N$ that share an ancestor. Let $u$ be an ancestor of both $v_1$ and $v_2$ in $N$. Suppose for contradiction that $L_{v_1} \cap L_{v_2} \neq \emptyset$, and let $x \in L_{v_1} \cap L_{v_2}$. Since $v_1$ and $v_2$ are descendant of $u$, we have $x \in L_u$. Now, let $P^1$ be a path from $u$ to $x$ containing $v_1$ and let $P^2$ be a path from $u$ to $x$ containing $v_2$. Since $v_1$ and $v_2$ are incomparable, $v_1$ is not a vertex of $P^2$. Hence, $P^1$ and $P^2$ are distinct, which is impossible in view of Lemma~\ref{lm:uniqp}.
\end{proof}

Next, we remark that} all leaves of $N$ are descendant of at least one root of $N$. In particular, if $N$ is a phylogenetic network, all leaves are descendant of the single root $r$ of $N$, that is, $L_r=L(N)$. In particular, in a phylogenetic network, any two leaves share an ancestor. However, this is not the case, in general, if $|R(N)| \geq 2$. For example the leaves $1$ and $5$ in the arboreal network $N$ depicted in Figure~\ref{fig-ntw}~(ii) do not share an ancestor. To keep track of the pairs of leaves that share an ancestor in $N$, the \emph{shared ancestry graph} $\mathcal A(N)$ of $N$, was introduced in \cite{HMS24}, and is defined as follows: The vertex set of $\mathcal A(N)$ is $L(N)$, and two vertices $x,y \in L(N)$ are joined by an edge in $\mathcal A(N)$ if and only if $x$ and $y$ share an ancestor in $N$ (see Figure~\ref{fig-shA} for the shared ancestry graphs of the networks depicted in Figure~\ref{fig-ntw}).

It was shown in \cite{HMS24} that, for all connected, undirected graph $G$, there exists a network $N$ with leaf set $V(G)$ such that $G=\mathcal A(N)$. Moreover, we have:

\begin{theorem}[\cite{HMS24}, Theorem~6.4]\label{th:ptol}
Let $G=(X,E)$ be an undirected graph. There exists an arboreal network $N$ such that $G=\mathcal A(N)$ if and only if $G$ is connected and Ptolemaic.
\end{theorem}

Note that a similar result was already established in \cite{UU09}, under a slightly different form. The main reason for referring to the result from \cite{HMS24} through this contribution, rather than to the one in \cite{UU09} is that the former is stated in terms of arboreal networks, which is the structure we are interested in here.

In this contribution, we focus on \emph{labelled networks}, that is, a pair $(N,t)$ such that $N$ is a network, and $t$ is a map from $V^*(N)$ to the set $\{0,1\}$. If $N$ is arboreal then $(N,t)$ induces an undirected graph $\mathcal C(N,t)$ with vertex set $L(N)$, such that two distinct elements $x,y \in L(N)$ are joined by an edge in $\mathcal C(N,t)$ if and only if the following two properties hold:
\begin{itemize}
\item[(1)] $x$ and $y$ share an ancestor in $N$.
\item[(2)] $t(\lca_N(x,y))=1$.
\end{itemize}

Note that by Lemma~\ref{lm:lcas}, $\lca_N(x,y) \in V^*(N)$, so $t(\lca_N(x,y))$ is well defined. Moreover, $\mathcal C(N,t)$ is uniquely determined by $(N,t)$.

In view of this, we say that an undirected graph $G=(X,E)$ is \emph{arboreal-explainable} if there exists a labelled arboreal network $(N,t)$ such that $G=\mathcal C(N,t)$. In that case, we say that $(N,t)$ \emph{explains} $G$.

Recall that labelled trees are special cases of labelled arboreal networks, since a labelled tree is a labelled arboreal network with one root. As mentioned in the introduction, an undirected graph $G$ is a cograph if and only if there exists a labelled phylogenetic tree $(T,t)$ explaining $G$. As a consequence, we have:

\begin{observation}\label{lm:cog}
Cographs are arboreal-explainable.
\end{observation}

Moreover, if $(N,t)$ is a labelled arboreal network such that $t(v)=1$ for all $v \in V^*(N)$, then one can easily verify that $\mathcal A(N)$ and $\mathcal C(N,t)$ are isomorphic. In particular, by Theorem~\ref{th:ptol} we have:

\begin{observation}\label{lm:ptol}
Connected Ptolemaic graphs are arboreal-explainable.
\end{observation}

We conclude this section with a technical result, which will be used in subsequent sections of the paper.

\begin{lemma}\label{lm:indcog}
Let $(N,t)$ be a labelled arboreal network, and let $G=\mathcal C(N,t)$. We have:
\begin{itemize}
\item[(i)] For all vertices $v$ of $N$, $G[L_v]$ is a cograph.
\item[(ii)] If $Y \subseteq V(G)$ is such that $G[Y]$ is not a cograph, there exists $x,y \in Y$ and $h \in H(N)$ such that $x \notin L_h$, $y \in L_h$, and $\{x,y\}$ is an edge of $G$.
\end{itemize}
\end{lemma}

\begin{proof}
(i) Let $v$ be a vertex of $N$. Without loss of generality, we may assume that $v$ has outdegree at least $2$ in $N$. Indeed, if $v$ has outdegree $0$, then $|L_v|=1$ and $G[L_v]$ is a cograph. If $v$ has outdegree $1$ in $N$, then $L_v=L_{v'}$ holds for $v'$ the unique child of $v$, so we can choose $v'=v$ and repeat the process if $v'$ has outdegree $1$. Let $N_v$ be the graph obtained from $N$ by removing all vertices of $N$ that are not descendant of $v$ and their incident arcs, and suppressing resulting vertices with indegree and outdegree $1$. By construction, $N_v$ is a network with $R(N_v)=\{v\}$, and $L(N_v)=L_v$. Moreover, since $N$ is arboreal and $N_v$ is a subgraph of $N$, $N_v$ is arboreal. In particular, $N_v$ is a phylogenetic tree.

Now, consider the restriction $t_v$ of $t$ to $V^*(N_v)$. Note that $t_v$ is well defined, since by construction, all vertices of $V^*(N_v)$ are vertices of $V^*(N)$. To see that $G[L_v]$ is a cograph, it suffices to show that $(N_v,t_v)$ explains $G[L_v]$. So, let $x,y \in L_v$ distinct. Note that $x$ and $y$ share an ancestor in $N$, since they are both descendant of $v$. Moreover, we have $\lca_{N_v}(x,y)=\lca_{N}(x,y)$ by construction of $N_v$, and $t_v(\lca_{N_v}(x,y))=t(\lca_{N}(x,y))$ by definition of $t_v$. Since $(N,t)$ explains $G$, it follows that $x$ and $y$ are joined by an edge in $G[L_v]$ if and only if $t_v(\lca_{N_v}(x,y))=1$. Hence, $(N_v,t_v)$ explains $G[L_v]$, so $G[L_v]$ is a cograph.

(ii) Let $Y \subseteq V(G)$ is such that $G[Y]$ is not a cograph. Then, there exists $x,y,z,u \in Y$ distinct such that $G[\{x,y,z,u\}]$ is a $P_4$. Up to a permutation, we may assume that $G[\{x,y,z,u\}]$ has edge set\\ $\{\{x,y\},\{y,z\},\{z,u\}\}$. Since $(N,t)$ explains $G$, it follows that $x$ and $y$ share an ancestor in $N$, $y$ and $z$ share an ancestor in $N$, and $z$ and $u$ share an ancestor in $N$. Put $v_1=\lca_N(x,y)$, $v_2=\lca(y,z)$, and $v_3=\lca(z,u)$. Note that $G[\{x,y,z,u\}]$ is not a cograph, so, in view of (i), there exists no $v \in V(N)$ such that $\{x,y,z,u\} \subseteq L_v$. In particular, \new{$v_1 \neq v_3$ must hold}. We next claim that if $v_1$ and $v_2$ share an ancestor in $N$, then $v_1=v_2$ must hold. 

To see the claim, suppose that $v_1$ and $v_2$ share an ancestor in $N$. Note that $y \in L_{v_1} \cap L_{v_2}$. \new{By Lemma~\ref{lm:emptyint}, $v_1$ and $v_2$ cannot be incomparable, so one must be an ancestor of the other}. If $v_1$ is a proper ancestor of $v_2$ in $N$, then \new{$v_1$ has a child $w$ that is an ancestor of $v_2$. In particular, $y,z \in L_w$ and, since $v_1=\lca_N(x,y)$, $z \notin L_w$. By Lemma~\ref{lm:uniqp}, $v_1$ does not have any further child $w'$ with $x \in L_{w'}$, so $v_1=\lca_N(x,z)$ follows.} However, $(N,t)$ explains $G$, and $\{x,y\}$ is an edge of $G$ while $\{x,z\}$ is not, so we have $t(\lca_N(x,y))=1$ and $t(\lca_N(x,z))=0$. \new{Since $\lca_N(x,y)=\lca_N(x,z)=v_1$, this is impossible}. If otherwise, $v_2$ is a proper ancestor of $v_1$ in $N$, then \new{by similar arguments,} we have $v_2=\lca_N(x,z)$. Again, $(N,t)$ explains $G$, and $\{y,z\}$ is an edge of $G$ while $\{x,z\}$ is not, so we have $t(\lca_N(y,z))=1$ and $t(\lca_N(x,z))=0$. \new{Since $\lca_N(x,z)=\lca_N(y,z)=v_2$, this is impossible}. To summarize, $v_1$ and $v_2$ are not incomparable in $N$, and neither vertex is a proper ancestor of the other, so $v_1=v_2$ must hold, which proves the claim.

By symmetry, the same holds for $v_2$ and $v_3$, that is, if $v_2$ and $v_3$ share an ancestor in $N$, then $v_2=v_3$ must hold. Since $v_1 \neq v_3$, it follows that $v_2$ cannot share an ancestor with both $v_1$ and $v_3$. Without loss of generality, we may assume that $v_1$ and $v_2$ do not share an ancestor in $N$. Since both $v_1$ and $v_2$ are ancestors of $y$, there exists a path $P^1$ from $v_1$ to $y$ in $N$, and a path $P^2$ from $v_2$ to $y$ in $N$. Since $v_1$ and $v_2$ do not share an ancestor in $N$, we have in particular that $v_1 \notin V(P^2)$ and $v_2 \notin V(P^1)$. Moreover, we have $y \in V(P^1) \cap V(P^2)$. In particular,  $V(P^1) \cap V(P^2) \neq \emptyset$. Now, let $h$ be the (necessarily unique) vertex such that $h \in V(P^1) \cap V(P^2)$, and no strict ancestor of $h$ enjoys that property. The choice of $h$, together with that fact that $h$ is distinct from $v_1$ and $v_2$, implies that $h$ has one parent in $V(P^1)$ and one parent in $V(P^2)$, and these parents are distinct. Hence, $h \in H(N)$. Clearly, $y \in L_h$. Moreover, since $v_1=\lca_N(x,y)$, no vertex of $P^1$ distinct from $v_1$ is an ancestor of $x$. In particular, $x \notin L_h$. Since $\{x,y\}$ is an edge of $G$, this concludes the proof.
\end{proof}

\section{Connectedness in arboreal-explainable graphs}\label{sec-con}

In Observation~\ref{lm:ptol}, we stated that connected Ptolemaic graphs are arboreal-explainable. In this section, we show that an undirected graph $G$ is arboreal-explainable if and only if all connected components of $G$ are arboreal-explainable (Proposition~\ref{pr:red}). In particular, the connectedness condition of Observation~\ref{lm:ptol} can be dropped. We start with a useful lemma.

\begin{lemma}\label{lm:ccc}
Let $G_1$ and $G_2$ be two vertex-disjoint undirected graphs, and let $G$ be the disjoint union of $G_1$ and $G_2$. If there exists a labelled arboreal network $(N_1,t_1)$ explaining $G_1$, a labelled arboreal network $(N_2,t_2)$ explaining $G_2$, then there exists a labelled arboreal network $(N,t)$ explaining $G$ such that $|R(N)|=|R(N_1)|+|R(N_2)|-1$.
\end{lemma}

\begin{proof}
Consider the network $N$ obtained from $N_1$ and $N_2$ by adding a new vertex $\rho$ with children $\rho_1$ and $\rho_2$, where $\rho_1$ and $\rho_2$ are two roots of $N_1$ and $N_2$, respectively. \new{Clearly, $N$ is arboreal, and we have $L(N)=L(N_1) \cup L(N_2)=V(G_1) \cup V(G_2)=V(G)$.} We define a labelling $t$ of $V^*(N)$ by putting $t(v)=t_1(v)$ if $v \in V^*(N_1)$, $t(v)=t_2(v)$ if $v \in V^*(N_2)$, and $t(v)=0$ otherwise, that is, if $v=\rho$. By construction, $N$ has exactly $|R(N_1)|+|R(N_2)|-1$ roots: all roots of $N_1$ distinct from $\rho_1$, all roots of $N_2$ distinct from $\rho_2$, and $\rho$. It remains to show that $(N,t)$ explains $G$.

To do this, consider two vertices $x$ and $y$ of $G$. If $x$ and $y$ belong to $V(G_i)$ for some $i \in \{1,2\}$, then $x$ and $y$ are leaves of $N_i$, and $x$ and $y$ share an ancestor in $N$ if and only if they share an ancestor in $N_i$. If this is the case, then we have $\lca_N(x,y)=\lca_{N_i}(x,y)$ by construction of $N$, and $t(\lca_N(x,y))=t_i(\lca_{N_i}(x,y))$ by definition of $t$. Since $(N_i,t_i)$ explains $G_i$, it follows that $\{x,y\}$ is an edge of $G$ if and only if $x$ and $y$ share an ancestor in $N$ and $t(\lca_N(x,y))=1$.

If otherwise $x \in V(G_1)$ and $y \in V(G_2)$, then \new{since $G$ is the disjoint union of $G_1$ and $G_2$,} $x$ and $y$ are not joined by an edge in $G$. Thus, we need to show that, if $x$ and $y$ share an ancestor in $N$, then $t(\lca_N(x,y))=0$. So, suppose that $x$ and $y$ share an ancestor in $N$. Recall that $x$ is a vertex of $N_1$, and $y$ is a vertex of $N_2$. By construction, there is only one vertex in $N$ that has descendants in both $N_1$ and $N_2$, and that vertex is $\rho$. Hence, $\lca_N(x,y)=\rho$ must hold. Since $t(\rho)=0$ by definition of $t$, we have $t(\lca_N(x,y))=0$ as desired. This concludes the proof that $(N,t)$ explains $G$.
\end{proof}

As a consequence, we have:

\begin{proposition}\label{pr:red}
Let $G=(X,E)$ be an undirected graph. Then $G$ is arboreal-explainable if and only if all connected components of $G$ are arboreal-explainable.
\end{proposition}

\begin{proof}
Suppose first that $G$ is arboreal-explainable, and let $(N,t)$ be a labelled arboreal network explaining $G$. Let $H$ be a connected component of $G$, and put $Y=V(H)$. We show that $H$ is arboreal-explainable. Consider the network $N'$ obtained from $N$ by successively:
\begin{itemize}
\item[(1)] Removing all vertices $v \in V(N)$ such that $L_v \cap Y=\emptyset$ and their incident arcs.
\item[(2)] Removing resulting vertices of indegree $0$ and outdegree $1$ and their incident arc.
\item[(3)] Suppressing resulting vertices of indegree and outdegree $1$.
\end{itemize}

By construction, we have $L(N')=Y$, and $V^*(N') \subseteq V^*(N)$. Hence, we can define the labelling $t'$ of $V^*(N')$ as the restriction of $t$ to $V^*(N')$. We first show that $N'$ is connected. To see this, suppose for contradiction that this is not the case, that is, $N$ has two distinct connected components $N_1$ and $N_2$. Since $H$ is connected, and $L(N_1) \cup L(N_2) \subseteq L(N')=Y$, there exists $x \in L(N_1)$ and $y \in L(N_2)$ such that $\{x,y\}$ is an edge of $H$. Since $H$ is a subgraph of $G$, and $(N,t)$ explains $G$, this means that $x$ and $y$ share an ancestor $u$ in $N$. In particular, $u$ is a vertex of $N'$ and is an ancestor of both $x$ and $y$ in $N'$, a contradiction since $x$ and $y$ belong to distinct connected components of $N'$. Hence, $N'$ is connected. Moreover, since $N'$ is a subgraph of $N$, and $N$ is arboreal, then the underlying undirected graph of $N'$ \new{does not contain a cycle}. Hence, $N'$ is an arboreal network with leaf set $Y$. It remains to show that $(N',t')$ explains $H$.

So, let $x,y \in Y$ distinct. If $x$ and $y$ do not share an ancestor in $N$, then $x$ and $y$ do not share an ancestor in $N'$. Otherwise, by construction of $N'$, we have $\lca_{N'}(x,y)=\lca_N(x,y)$, and $t'(\lca_{N'}(x,y))=t(\lca_N(x,y))$ holds by definition of $t'$. Since $(N,t)$ explains $G$, and $H$ is an induced subgraph of $G$, it follows that $x$ and $y$ are joined by an edge in $H$ if and only if $x$ and $y$ share an ancestor in $N'$, and $t'(\lca_{N'}(x,y))=1$. Hence, $(N',t')$ explains $H$, which concludes the proof that $H$ is arboreal-explainable.

Conversely, suppose that all connected components of $G$ are arboreal-explainable. \new{If $G$ is connected, then there is nothing to show. Otherwise}, repeated applications of Lemma~\ref{lm:ccc} imply that $G$ is arboreal-explainable.
\end{proof}

Interestingly, there is a simple charaterization of labelled arboreal networks $(N,t)$ for which $\mathcal C(N,t)$ is connected.

\begin{proposition}
Let $(N,t)$ be a labelled arboreal network with $|L(N)| \geq 2$, and let $G=\mathcal C(N,t)$. Then $G$ is connected if and only if $t(r)=1$ holds for all $r \in R(N)$.
\end{proposition}

\begin{proof}
Let $X=L(N)$.

Suppose first that there exists a root $r \in R(N)$ such that $t(r)=0$, and let $k \geq 2$ be the degree of $r$ in $N$. Let $N'$ be the directed graph obtained from $N$ by removing $r$ and its incidents arcs. Then $N'$ has $k \geq 2$ connected components, each containing one of the children of $r$ in $N$. Now, let $x,y \in L(N)$ be such that $x$ and $y$ belong to two distinct connected components of $N'$. If $x$ and $y$ do not share an ancestor in $N$, then since $(N,t)$ explains $G$, $x$ and $y$ are not joined by an edge in $G$. If otherwise, $x$ and $y$ share an ancestor in $N$, then $\lca_N(x,y)=r$ must hold. In that case, since $t(r)=0$ and $(N,t)$ explains $G$, $x$ and $y$ are not joined by an edge in $G$. It follows that, for any two connected components $N_1$ and $N_2$ of $N'$ there is no edge in $G$ joining a vertex of $V(N_1) \cap X$ with a vertex of $V(N_2) \cap X$. Hence, $G$ has at least $k \geq 2$ connected components, so $G$ is disconnected.

Conversely, suppose that $t(r)=1$ holds for all $r \in R(N)$. Let $x,y \in V(G)$ distinct. We show that $x$ and $y$ are connected in $G$. Suppose first that $x$ and $y$ share an ancestor in $N$, and let $v=\lca_N(x,y)$. If $t(v)=1$, then since $N$ explains $G$, $\{x,y\}$ is an edge of $G$. If otherwise, $t(v)=0$, then by assumption on $N$, $v$ is not a root of $N$. Now, let $r$ be a root of $N$ that is an ancestor of $v$. \new{Let $u$ be a child of $r$ that is an ancestor of $v$, and let $u'$ be a further child of $r$ (which must exist since $r$ has outdegree at least two). By Lemma~\ref{lm:uniqp}, we have $L_u \cap L_{u'}=\emptyset$. Moreover, we have $x,y \in L_v \subseteq L_u$. Thus, for} all $z \in L_{u'}$, we have $\lca_N(x,z)=\lca_N(y,z)=r$. In particular, $t(\lca_N(x,z))=t(\lca_N(y,z))=t(r)=1$, so $\{x,z\}$ and $\{y,z\}$ are edges of $G$, so $x$ and $y$ are connected in $G$.

Suppose now that $x$ and $y$ do not share an ancestor in $G$. By Theorem~\ref{th:ptol}, the shared ancestry graph $\mathcal A(N)$ of $N$ is connected. Hence, there exists a sequence $x_1, \ldots, x_k$, $k \geq 2$ of elements of $G$ such that $x_1=x$, $x_k=y$, and for all $i \in \{1, k-1\}$, $x_i$ and $x_{i+1}$ share an ancestor in $N$. In view of the above, for all $i \in \{1, k-1\}$, $x_i$ and $x_{i+1}$ are connected in $G$. Hence $x$ and $y$ are connected in $G$, which concludes the proof that $G$ is connected.
\end{proof}

\section{Characterization}\label{sec-char}

In this section, we present the main result of the paper (Theorem~\ref{th:dh}), namely, that an undirected graph $G$ is arboreal-explainable if and only if $G$ is distance-hereditary.

We start with showing that distance-hereditary graphs are arboreal-explainable. More specifically, we show that the recursive definition of a distance-hereditary graph $G$ given in Section~\ref{sec-ug} can be used to build a labelled arboreal network $(N,t)$ explaining $G$. We begin with some definitions.

For $G=(X,E)$ a graph and $x_1, \ldots, x_k$, $k =|X|$ an ordering of the elements of $X$, we define $G_i=G[\{x_1, \ldots, x_i\}]$ for all $i \in \{1, \ldots, k\}$. \new{In particular,} $G_1$ is the singleton graph $\{x_1\}$, and $G_k=G$. Recall from Section~\ref{sec-ug} that $G$ is distance-hereditary if and only if there exists an ordering $x_1, \ldots, x_k$ of the elements of $X$ such that, for all $i \in \{1, \ldots, k-1\}$, $G_{i+1}$ can be obtained from $G_i$ by adding $x_{i+1}$ as (R1) a pendant-vertex, (R2) a false-twin or (R3) a true-twin to some vertex $x \in \{x_1, \ldots, x_i\}=V(G_i)$. Following \cite{BM86}, we call such an ordering a \emph{sequence of one-vertex extensions}. It was shown in \cite{HM:90}, that, if $G$ is distance-hereditary, such an sequence can be found in linear time (see also \cite{DHP01} for an updated version of the algorithm).

\new{Given a distance-hereditary graph $G$, the following algorithm uses this concept of a sequence of one-vertex extensions to build a labelled arboreal network explaining $G$.}

\begin{algorithm}[h]
  \small 
  \caption{\texttt{From distance-hereditary graph to labelled arboreal network}}
  \label{alg:dh}
  \begin{algorithmic}[1]
    \Require  A distance-hereditary graph $G=(X,E)$ with $|X| \geq 2$.
    \Ensure   A labelled aboreal network $(N,t)$ explaining $G$.

   \State Put $k=|X|$.
    \State Find a sequence $x_1, \ldots, x_k$ of one vertex-extensions for $G$. \label{l:funky}
   \State Initialize $N$ as a directed graph with vertex set $\{v_0,x_1,x_2\}$ and arc set $\{(v_0,x_1),(v_0,x_2)\}$.\label{l:initN}
    \If{$\{x_1,x_2\}$ is an edge of $G$.} \label{l:12e}
    	\State Put $t(v_0)=1$. \label{l:t01}
\Else \label{l:12ne}
    	\State Put $t(v_0)=0$. \label{l:t00}
   \EndIf
   \If {$k \geq 3$} \label{l:if3}
    \For{$i$ from $3$ to $k$} \label{l:for}    
    	\If{$G_i$ is obtained from $G_{i-1}$ by adding $x_i$ as a false- or true-twin to some vertex $x \in V(G_{i-1})$.} \label{l:R23}
    	\State Subdivide the incoming arc $e$ of $x$ in $N$ by introducing a new vertex $v$. \label{l:sub2}
    		\State Add to $N$ the arc $(v,x_i)$. \label{l:add2}
    		\If{$\{x,x_i\}$ is an edge of $G$.} \label{l:ie} \Comment{$x$ and $x_i$ are true-twins}
    			\State Put $t(v)=1$. \label{l:ti1}
			\Else \label{l:nie} \Comment{$x$ and $x_i$ are false-twins}
    			\State Put $t(v)=0$. \label{l:ti0}
  			\EndIf
	\Else \Comment{$G_i$ is obtained from $G_{i-1}$ by adding $x_i$ as a pendant-vertex to some vertex $x \in V(G_{i-1})$.} \label{l:R1}
			\State Subdivide the incoming arc $e$ of $x$ in $N$ by introducing a new vertex $v$. \label{l:sub1}
    		\State Add to $N$ a new vertex $u$, and the arcs $(u,v)$ and $(u,x_i)$. \label{l:add1}
    		\State Put $t(u)=1$. \label{l:tp}
    \EndIf
	\EndFor \label{l:forend}
	\EndIf
	\State \Return $(N,t)$.
  \end{algorithmic}
\end{algorithm}

We next show the correctness of Algorithm~\ref{alg:dh}.

\begin{proposition}\label{pr:algo}
For $G=(X,E)$ a distance-hereditary graph with $|X| \geq 2$, Algorithm~\ref{alg:dh} correctly computes a labelled arboreal network $(N,t)$ explaining $G$.
\end{proposition}

\begin{proof}
Let $(N_2, t_2)$ be the labelled network built at Lines~\ref{l:initN} to \ref{l:t00}. Note that if $k=2$, the loop starting at Line~\ref{l:for} is not visited, and we have $(N,t)=(N_2,t_2)$. If $k \geq 3$, then for $i \in \{3, \ldots, k\}$, let $(N_i,t_i)$ be the network obtained after the $(i-2)$th iteration of the loop at Line~\ref{l:for} (that is, after the leaf $x_i$ has been added to the network). Since the loop at Line~\ref{l:for}, ends after $k-2$ iterations, we have $(N,t)=(N_{k},t_{k})$. We next show that, for all $i \in \{2, \ldots, k\}$, $(N_i,t_i)$ is a labelled arboreal network explaining $G_i$.

\new{Consider first the case $i=2$. The network} $N_2$ is the phylogenetic tree with three vertices, $v_0$, $x_1$, and $x_2$, and two arcs $(v_0,x_1)$ and $(v_0,x_2)$ (Line~\ref{l:initN}). In particular, $N_2$ is an arboreal network, and $L(N_2)=\{x_1,x_2\}=V(G_2)$. Moreover we have $t_2(v_0)=1$ if $\{x_1,x_2\}$ is an edge of $G$, and $t_2(v_0)=0$ otherwise (Lines~\ref{l:12e} to \ref{l:t00}). Since $v_0=\lca_{N_2}(x_1,x_2)$, it follows that $(N_2,t_2)$ explains $G_2$.

\new{If $k=2$, then we are done. If otherwise, $k \geq 3$, we show that $(N_i,t_i)$ is a labelled arboreal network explaining $G_i$ for all $i \in \{2, \ldots, k\}$ by induction on $i$. The case $i=2$ being established, we now assume that $i \geq 3$, and} that $(N_{i-1},t_{i-1})$ is a labelled arboreal network explaining $G_{i-1}$. \new{Recall that $x_1, \ldots, x_k$ is a sequence of one vertex-extensions for $G$. In particular, $G_i$ is obtained from $G_{i-1}$ by adding the vertex $x_i$ either as a pendant-vertex, or as a false- or true-twin to some vertex $x$ of $G_{i-1}$. We investigate these two cases in turn.}

Suppose first that $G_i$ is obtained from $G_{i-1}$ by adding $x_i$ as a pendant-vertex to some vertex $x \in V(G_{i-1})$. Then $N_i$ is obtained from $N_{i-1}$ by subdividing the incoming arc $e$ of $x$ in $N_{i-1}$ with the introduction of a new vertex $v$ (Line~\ref{l:sub1}), and adding a new vertex $u$, and the arcs $(u,v)$ and $(u,x_i)$ (Line~\ref{l:add1}).

Since $N_{i-1}$ is an arboreal network by assumption, $N_i$ is an arboreal network. Indeed, the operation described above does not disconnect the graph, and does not create a cycle in the underlying, undirected graph. Moreover, all vertices of $N_i$ have the same indegree and outdegree as in $N_{i-1}$, and the new vertices $v,u$, and $x_i$ have indegree $2$ and outdegree $1$, indegree $0$ and outdegree $2$, and indegree $1$ and outdegree $0$, respectively. In particular, we have $V^*(N_i)=V^*(N_{i-1}) \cup \{u\}$, and $L(N_i)=L(N_{i-1}) \cup \{x_i\}=V(G_i)$. Since \new{$t_i=t_{i-1}$ on  $V^*(N_{i-1})$, and $t_i(u)$ is defined at Line~\ref{l:tp}}, the domain of $t_i$ is $V^*(N_{i-1}) \cup \{u\} \new{=V^*(N_i)}$, so $(N_i,t_i)$ is a labelled arboreal network.

To see that $(N_i,t_i)$ explains $G_i$, let $y,z$ be two distinct elements of $V(G_i)$. If neither $y=x_i$ nor $z=x_i$ holds, then we have $y,z\in L(N_{i-1})$. In this case,  $y$ and $z$ share an ancestor in $N_i$ if and only if $y$ and $z$ share an ancestor in $N_{i-1}$. Moreover, if this holds, $\lca_{N_i}(y,z)=\lca_{N_{i-1}}(y,z)$, and $t_i(\lca_{N_i}(y,z))=t_{i-1}(\lca_{N_{i-1}}(y,z))$, since the value of $t$ on the elements of $V^*(N_{i-1})$ is not modified. Since $(N_{i-1},t_{i-1})$ explains $G_{i-1}$ by our induction hypothesis, \new{and $\{y,z\}$ is an edge of $G_i$ if and only if $\{y,z\}$ is an edge of $G_{i-1}$,} it follows that $\{y,z\}$ is an edge of $G_i$ if and only if $y$ and $z$ share an ancestor in $N_i$ and $t_i(\lca_{N_i}(y,z))=1$. \new{Suppose now that one of} $y=x_i$ or $z=x_i$ holds, say $y=x_i$. \new{By construction,} $y$ and $z$ share an ancestor in $N$ if and only if $z=x$. Moreover, $G_i$ is obtained from $G_{i-1}$ by adding $x_i$ as a pendant edge to $x$. In particular, $\{x_i,z\}$ is an edge of $G_i$ if and only if $z=x$. By construction, we have $\lca_{N_i}(x_i,x)=u$, and we have $t_i(u)=1$ (Line~\ref{l:tp}). \new{Taken together, these observations imply that $\{x_i,z\}$ is an edge of $G_i$ if and only if $x_i$ and $z$ share an ancestor in $N$ and $t(\lca_N(x_i,z))=1$. This} concludes the proof that $(N_i,t_i)$ explains $G_i$.

Suppose now that $G_i$ is obtained from $G_{i-1}$ by adding $x_i$ as a false- or true-twin to some vertex $x \in V(G_{i-1})$. Then $N_i$ is obtained from $N_{i-1}$ by subdividing the incoming arc $e$ of $x$ in $N_i$ with the introduction of a new vertex $v$ (Line~\ref{l:sub2}), and adding the arc $(v,x_i)$ (Line~\ref{l:add2}).

Using similar arguments as in the previous case, one can see that $N_i$ is an arboreal network, that $V^*(N_i)=V^*(N_{i-1}) \cup \{v\}$, and that $L(N_i)=L(N_{i-1}) \cup \{x_i\}=V(G_i)$. Since \new{$t_i=t_{i-1}$ on  $V^*(N_{i-1})$, and $t_i(v)$ is defined at Line~\ref{l:ti1} or \ref{l:ti0}}, the domain of $t_i$ is $V^*(N_{i-1}) \cup \{v\}\new{=V^*(N_i)}$, so $(N_i,t_i)$ is a labelled arboreal network.

To see that $(N_i,t_i)$ explains $G_i$, let $y,z$ be two distinct elements of $V(G_i)$. If neither $y=x_i$ nor $z=x_i$ holds, then similar arguments as in the previous case show that $\{y,z\}$ is an edge of $G_i$ if and only if $y$ and $z$ share an ancestor in $N_i$ and $t_i(\lca_{N_i}(y,z))=1$. \new{Suppose now that one of} $y=x_i$ or $z=x_i$ holds, say $y=x_i$. \new{By construction,}  $y$ and $z$ share an ancestor in $N_i$ if and only if $z$ and $x$ share an ancestor in $N_{i-1}$. We next distinguish between three cases: (a) $z=x$ (b) $z$ and $x$ do not share an ancestor in $N_{i-1}$ and (c) \new{$z \neq x$ and} $z$ and $x$ share an ancestor in $N_{i-1}$.

\begin{itemize}
\item[(a)] $z=x$. In this case, we have $\lca_{N_i}(x_i,z)=v$, and we have $t_i(v)=1$ if and only if $\{x_i,z\}$ is an edge of $G_i$ (Lines~\ref{l:ie} to \ref{l:ti0}).
\item[(b)] $z$ and $x$ do not share an ancestor in $N_{i-1}$. \new{In particular, $z \neq x$.} In this case, $z$ and $x_i$ do not share an ancestor in $N_i$, so we need to show that $\{x_i,z\}$ is not an edge of $G_i$. Since $z$ and $x$ do not share an ancestor in $N_{i-1}$, and $(N_{i-1},t_{i-1})$ explains $G_{i-1}$, it follows that $\{x,z\}$ is not an edge in $G_{i-1}$. Moreover, $G_i$ is obtained from $G_{i-1}$ by adding $x_i$ as a false- or true-twin of $x$, so \new{the fact that $\{x,z\}$ is not an edge in $G_{i-1}$ implies that} $\{x_i,z\}$ is not an edge of $G_i$, \new{as desired}.
\item[(c)] \new{$z \neq x$ and} $z$ and $x$ share an ancestor in $N_{i-1}$. In this case, $z$ and $x_i$ share an ancestor in $N_i$. \new{Hence,} we need to show that $\{x_i,z\}$ is an edge of $G_i$ if and only if $t_i(\lca_{N_i}(x_i,z))=1$. \new{Since  $G_i$ is obtained from $G_{i-1}$ by adding $x_i$ as a false- or true-twin of $x$, $\{x_i,z\}$ is an edge of $G_i$ if and only if $\{x,z\}$ is an edge of $G_{i-1}$. Suppose first that $\{x_i,z\}$ is an edge of $G_i$. Then, $\{x,z\}$ is an edge of $G_{i-1}$, and since $(N_{i-1},t_{i-1})$ explains $G_{i-1}$ by our induction hypothesis, we have $t_{i-1}(\lca_{N_{i-1}}(x,z))=1$. Moreover, we have $\lca_{N_{i-1}}(x,z)=\lca_{N_{i}}(x,z)=\lca_{N_{i}}(x_i,z)$ by construction of $N$, and $t_i(\lca_{N_{i}}(x,z))=t_{i-1}(\lca_{N_{i-1}}(x,z))$ by definition of $t_i$, so $t_i(\lca_{N_{i}}(x_i,z))=1$ follows. Suppose now that $\{x,z\}$ is not an edge of $G_{i-1}$. Then similar arguments as in the previous case imply that $t_i(\lca_{N_{i}}(x_i,z))=0$.}
\end{itemize}

\new{To summarize, in all three cases, $\{x_i,z\}$ is an edge of $G_i$ if and only if $x_i$ and $z$ share an ancestor in $N$ and $t(\lca_N(x_i,z))=1$. This} concludes the proof that $(N_i,t_i)$ explains $G_i$. \new{Since $G_k=G$ and $(N,t)=(N_k,t_k)$, the conclusion follows.}
\end{proof}

\begin{figure}[h]
\begin{center}
\includegraphics[scale=0.7]{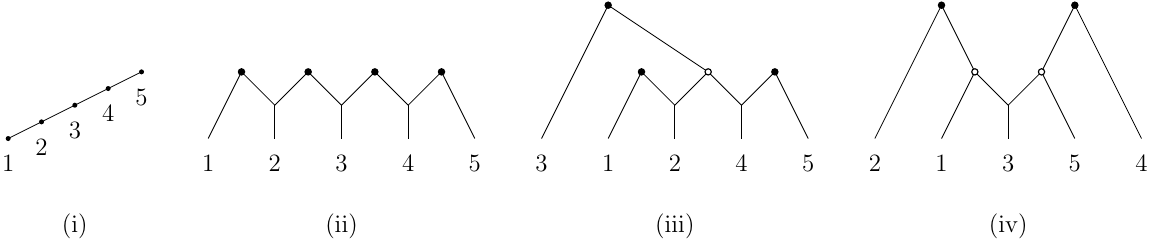}
\end{center}
\caption{(i) An undirected graph $G$ with vertex set $\{1,2,3,4,5\}$. (ii) to (iv) Three distinct binary labelled arboreal networks explaining $G$, where vertices with label $1$ are indicated with $\bullet$, and vertices with label $0$ are indicated with $\circ$.}
\label{fig-alg}
\end{figure}

Note that for $G$ a distance-hereditary graph, the labelled arboreal-network outputted by Algorithm~\ref{alg:dh} when given $G$ as an input may not be unique. The reason for this is that the construction depends on the sequence of one-vertex extensions computed at Line~\ref{l:funky}, which, in general, is not unique. For example, consider the undirected graph $G$ with vertex set $X=\{1,2,3,4,5\}$ depicted in Figure~\ref{fig-alg}~(i). Then, both $1,2,3,4,5$ and $3,2,4,1,5$ are sequences of one-vertex extensions for $G$. One can verify that Algorithm~\ref{alg:dh} returns the labelled arboreal network depicted in Figure~\ref{fig-alg}~(ii) in case the sequence computed at Line~\ref{l:funky} is $1,2,3,4,5$. If that sequence is $3,2,4,1,5$, then Algorithm~\ref{alg:dh} returns the labelled arboreal network depicted in Figure~\ref{fig-alg}~(iii).

Another interesting property of Algorithm~\ref{alg:dh} is that the outputted labelled network $(N,t)$ is always binary. Indeed, this is the case of the network $N$ built at Line~\ref{l:initN}. Then, at each iteration within the loop at Line~\ref{l:for}, new vertices are created, with either indegree $2$ and outdegree $1$ (vertex $v$, at Line~\ref{l:add1}), indegree $0$ and outdegree $2$ (vertex $u$, at Line~\ref{l:add1}), or indegree $1$ and outdegree $2$ (vertex $v$, at Line~\ref{l:add2}). At the same time, the indegree and outdegree of already present vertices are not modified.

Clearly, not all labelled arboreal networks representing a given distance-hereditary graph $G$ are binary, so not all such networks can be built using Algorithm~\ref{alg:dh}. More interestingly, there may exist binary arboreal networks $(N,t)$ explaining $G$ that can not be obtained using Algorithm~\ref{alg:dh}. This is the case, for example for the labelled arboreal network depicted in Figure~\ref{fig-alg}~(iv). Although that network is binary, and explains the undirected graph $G$ depicted in Figure~\ref{fig-alg}~(i), one can never obtain $(N,t)$ as an output of Algorithm~\ref{alg:dh}. To see this, it suffices to remark that, after each iteration of the loop at Line~\ref{l:for}, the current network contains either a vertex $v$ such that both children of $v$ are leaves (Lines~\ref{l:sub2} and \ref{l:add2}), or a vertex $u$ such that $u$ has indegree $0$, one child of $u$ is a leaf, and the other is a vertex of indegree $2$ (Lines~\ref{l:sub1} and \ref{l:add1}). Since Algorithm~\ref{alg:dh} finishes immediately after the last iteration of the loop, and the network depicted in Figure~\ref{fig-alg}~(iv) does not contain either of these two types of vertices, it follows that that network cannot be obtained using Algorithm~\ref{alg:dh}.

To complete the characterization of arboreal-explainable graphs, we first need a useful lemma.

\begin{lemma}\label{lm:rem}
Let $G=(X,E)$ be an arboreal-explainable graph, and let $(N,t)$ be a labelled arboreal network explaining $G$. \new{Let $x \in X$ be such that the parent $v$ of $X$ in $N$ has outdegree $2$ or more.} Then $G[X \setminus \{x\}]$ is arboreal-explainable.
\end{lemma}

\begin{proof}
Put $G_x=G[X \setminus \{x\}]$. Consider $N_x$ the network obtained from $N$ be removing $x$ and the arc $(v,x)$, removing any resulting vertex of indegree $0$ and outdegree $1$, and suppressing any resulting vertex of indegree and outdegree $1$. By construction, $N_x$ is an arboreal network with leaf set $L(N_x)=X \setminus \{x\}=V(G_x)$. Moreover, we have $V^*(N_x) \subseteq V^*(N)$. Hence, \new{the pair $(N_x,t_x)$, where $t_x$ is defined as} $t_x(u)=t(u)$ for all $u \in V^*(N_x)$, \new{is a labelled arboreal network. We now proceed to show that $(N_x,t_x)$ explains $G_x$.}

For $y,z \in X \setminus \{x\}$ distinct, $y$ and $z$ share an ancestor in $N_x$ if and only if $y$ and $z$ share an ancestor in $N$. Moreover, if that holds, we have $\lca_{N_x}(y,z)=\lca_N(y,z)$, and we have $t_x(\lca_{N_x}(y,z))=t(\lca_N(y,z))$ by definition \new{of $t_x$}. Since $(N,t)$ explains $G$, and $G_x$ is an induced subgraph of $G$, it follows that $\{y,z\}$ is an edge of $G$ if and only if $y$ and $z$ share an ancestor in $N_x$, and $t_x(\lca_{N_x}(y,z))=1$. Thus, $(N_x,t_x)$ explains $G$ \new{as desired}, so $G$ is arboreal-explainable.
\end{proof}

We are now ready to state the main result of this section.

\begin{theorem}\label{th:dh}
Let $G=(X,E)$ be an undirected graph. Then $G$ is arboreal-explainable if and only if $G$ is distance-hereditary.
\end{theorem}

\begin{proof}
Suppose first that $G$ is distance-hereditary. By Proposition~\ref{pr:algo}, Algorithm~\ref{alg:dh} applied to $G$ returns a labelled arboreal network $(N,t)$ explaining $G$, so $G$ is arboreal-explainable.

Conversely, suppose that $G$ is arboreal-explainable, and let $(N,t)$ be a labelled arboreal network explaining $G$. We show that $G$ is distance-hereditary by induction on $k=|X|$. If $k \leq 4$, then $G$ is distance-hereditary, since all graphs on $4$ or less vertices are distance-hereditary.

Suppose now that $k \geq 5$, and that all arboreal-explainable graphs $H$ with $|V(H)|<k$ are distance-hereditary. Consider the undirected graph $\overline N$ obtained from $N$ by ignoring the direction of the arcs, and suppressing resulting vertices of degree $2$. Since $N$ is arboreal, $\overline N$ is a tree, \new{and has leaf set $X$}. Moreover, $\overline N$ does not contain any vertex of degree $2$, so $\overline N$ contains a vertex $v$ that is adjacent to at least two leaves, and to at most one non-leaf vertex. Note that $V(\overline N) \subseteq V(N)$, so $v$ is a vertex of $N$. Moreover, \new{an element $x \in X$} adjacent to $v$ in $\overline N$ is either a child of $v$ in $N$, or there exists a root $r_x$ in $N$ with exactly two children, $x$ and $v$. Based on this observation, \new{and the fact that $v$ is adjacent to at least two leaves in $\overline N$, at least one of the following three cases holds: (a) at least two leaves adjacent to $v$ in $\overline N$ are children of $v$ in $N$, (b) at least two leaves adjacent to $v$ in $\overline N$ are not children of $v$ in $N$ and (c) neither (a) nor (b) holds. Note that in the latter case, $v$ is adjacent to exactly two leaves in $\overline N$, one of which is a child of $v$ in $N$ and the other is not. We now investigate these three cases in turns.}

\begin{itemize}
\item[(a)]: \new{at least two leaves adjacent to $v$ in $\overline N$ are children of $v$ in $N$. Let $x,y \in X$ be two such leaves}. In particular, $v \in V^*(N)$. Let $G_y=G[X \setminus \{y\}]$. By Lemma~\ref{lm:rem}, there exists a labelled arboreal network $(N',t')$ on $X \setminus \{y\}$ that explains $G_y$. Since $x$ and $y$ are both children of $v$ in $N$, for all $z \in X \setminus \{x,y\}$, $x$ and $z$ share an ancestor in $N$ if and only if $y$ and $z$ share an ancestor in $N$. Moreover, if that holds, then $\lca_N(x,z)=\lca_N(y,z)$. \new{Since $(N,t)$ explains $G$, this implies that} $\{x,z\}$ is an edge of $G$ if and only if $\{y,z\}$ is an edge of $G$. \new{Hence,} $G$ \new{can be} obtained from $G_y$ by adding $y$ as a true-twin to $x$ if $\{x,y\}$ is an edge of $G$, or as a false-twin to $x$ otherwise. Since \new{by our induction hypothesis,} $G_y$ is distance-hereditary, \new{it follows from Lemma~\ref{lm:indh} that} $G$ is distance-hereditary.

\item[(b)]: \new{at least two leaves adjacent to $v$ in $\overline N$ are not children of $v$ in $N$. Let $x,y \in X$ be two such leaves. Since $x$ and $y$ are not children of $v$,} there exists two \new{(necessarily distinct)} roots $r_x$ and $r_y$ of outdegree $2$, such that the children of $r_x$ are $x$ and $v$ and the children of $r_y$ are $v$ and $y$. Note that $r_x$ and $r_y$ are both vertices of $V^*(N)$. Let $G_y=G[X \setminus \{y\}]$ and $G_x=G[X \setminus \{x\}]$. By Lemma~\ref{lm:rem}, $G_y$ and $G_x$ are arboreal-explainable. Hence, by our induction hypothesis, $G_x$ and $G_y$ are distance-hereditary. Since $r_x$ is the only proper ancestor of $x$ in $N$, and $v$ is the only other child of $r_x$ in $N$, a vertex $z \in X \setminus \{x\}$ share an ancestor with $x$ if and only $z$ is a descendant of $v$. In this case, we have $\lca_N(x,z)=r_x$. By symmetry, a vertex $z \in X \setminus \{y\}$ share an ancestor with $y$ if and only $z$ is a descendant of $v$, and we have $\lca_N(y,z)=r_y$ in that case. We now distinguish between two subcases:

\begin{itemize}
\item[b1]: $t(r_x)=t(r_y)=1$. In view of the above observation, a vertex $z \in X \setminus \{x,y\}$ shares an ancestor with $x$ if and only if $z$ shares an ancestor with $y$, and if this holds, we have \new{$t(\lca_N(x,z))=t(r_x)=t(r_y)=t(\lca_N(y,z))$}. \new{Since $(N,t)$ explains $G$, this implies that $\{x,z\}$ is an edge of $G$ if and only if $\{y,z\}$ is an edge of $G$. In addition,} $x$ and $y$ do not share an ancestor in $N$, \new{so $\{x,y\}$ is not an edge of $G$. It follows from these observations} that $G$ \new{can be} obtained from $G_y$ by adding $y$ as a false-twin to $x$. Since $G_y$ is distance-hereditary, $G$ is distance-hereditary by Lemma~\ref{lm:indh}.

\item[b2]: at least one of $t(r_x)=0$ or $t(r_y)=0$, say $t(r_y)=0$. In this case, for all vertices $z \in X \setminus \{y\}$ such that $y$ and $z$ share an ancestor in $N$, we have $t(\lca_N(y,z))=0$. \new{Since $(N,t)$ explains $G$, this implies that} $y$ has degree $0$ in $G$, so $G$ is the disjoint union of $G_y$ and the \new{graph with vertex set $\{y\}$ and empty edge set}. Since $G_y$ is distance-hereditary, \new{it follows directly from the definition of distance-hereditary graphs that} $G$ is distance-hereditary.
\end{itemize}

\item[(c)]: neither (a) nor (b) holds. \new{As remarked above, this implies that $v$ is adjacent to exactly two leaves $x$ and $y$ in $\overline N$, such that (up to permutation),} $x$ is a child of $v$ in $N$, and there exists a root $r$ in $N$ with exactly two children, $y$ and $v$. Note that $r \in V^*(N)$. Let $G_y=G[X \setminus \{y\}]$ and $G_x=G[X \setminus \{x\}]$. By Lemma~\ref{lm:rem}, $G_y$ is arboreal-explainable. Hence, by our induction hypothesis, $G_y$ is distance-hereditary. Note first that a vertex $z \in X \setminus \{y\}$ shares an ancestor with $y$ in $N$ is and only if $z$ is a descendant of $v$ in $N$. If this is the case, then $\lca_N(z,y)=r$ must hold. We next distinguish between three subcases:

\begin{itemize}
\item[c1]: $t(r)=0$. In this case, for all vertices $z \in X \setminus \{y\}$ such that $y$ and $z$ share an ancestor in $N$, we have $t(\lca_N(y,z))=0$. In particular, $y$ has degree $0$ in $G$, so $G$ is the disjoint union of $G_y$ and the \new{graph with vertex set $\{y\}$ and empty edge set}. Since $G_y$ is distance-hereditary, \new{it follows directly from the definition of distance-hereditary graphs that} $G$ is distance-hereditary.

\item[c2]: $t(r)=1$ and $v$ has outdegree $1$ in $N$. In this case, $x$ is the unique child of $v$ in $N$. In particular, $x$ is the only element of $X \setminus \{y\}$ that shares an ancestor with $y$ in $N$. Since $\lca_N(x,y)=r$, $t(r)=1$, \new{and $(N,t)$ explained $G$,} it follows that $x$ is the only vertex adjacent to $y$ in $G$. Hence, $G$ \new{can be} obtained from $G_y$ by adding $y$ as a pendant-vertex to $x$. Since $G_y$ is distance-hereditary, $G$ is distance-hereditary by Lemma~\ref{lm:indh}.

\item[c3]: $t(r)=1$ and $v$ has outdegree $2$ \new{or more} in $N$. \new{Recall that $v$ is chosen such that it has exactly one non-leaf adjacent vertex in $\overline N$. Since $x$ and $y$ are the only leaves of $\overline N$ adjacent to $v$ in $\overline N$, $v$ has degree $3$ in $\overline N$. Next, we remark that all vertices of $\overline N$ have the same degree in $\overline N$ and in $N$. In particular, $v$ has degree $3$ in $N$, so $v$ has exactly one child other than $x$, and no parent other than $r$.}

\new{Since $v$ has outdegree $2$ in $N$}, $v \in V^*(N)$, so $G_x$ is arboreal-explainable in view of Lemma~\ref{lm:rem}. By our induction hypothesis, $G_x$ is distance-hereditary. We next remark that, since $\lca_N(x,y)=r$, $t(r)=1$, \new{and $(N,t)$ explains $G$,} $\{x,y\}$ is an edge of $G$. Moreover, \new{$r$ is the unique parent of $v$ in $N$, so} an element $z \in X \setminus \{x,y\}$ shares an ancestor with $x$ in $N$ if and only if \new{$z$ is a descendant of $v$. Since $r$ is the only ancestor of $y$ in $N$, the latter holds if and only if} $z$ shares an ancestor with $y$. \new{In summarry, $z$ shares an ancestor with $x$ if and only if $z$ shares an ancestor with $y$. In addition, for such an element $z$,} we have $\lca_N(y,z)=r$ and $\lca_N(x,z)=v$. If $t(v)=1$, \new{then we have $t(\lca_N(y,z))=t(r)=t(v)=t(\lca_N(x,z))$ for all $z \in X \setminus \{x,y\}$ such that $z$ shares an ancestor of $x$ in $N$. Since $(N,t)$ explains $G$, and in view of the above equivalence,} a vertex $z \in X \setminus \{x,y\}$ is adjacent to $y$ in $G$ if and only if $z$ is adjacent to $x$ in $G$. Hence, $G$ is obtained from $G_y$ by adding $y$ as a true-twin to $x$. Since $G_y$ is distance-hereditary, $G$ is distance-hereditary by Lemma~\ref{lm:indh}. If otherwise, $t(v)=0$, then for all vertices $z \in X \setminus \{x,y\}$ such that $x$ and $z$ share an ancestor in $N$, we have $t(\lca_N(y,z))\new{=t(v)}=0$. \new{Since $(N,t)$ explains $G$,} it follows thatt $y$ is the only vertex adjacent to $x$ in $G$. Hence, $G$ is obtained from $G_x$ by adding $x$ as a pendant-vertex to $y$. Since $G_x$ is distance-hereditary, $G$ is distance-hereditary by Lemma~\ref{lm:indh}.
\end{itemize}
\end{itemize}
\end{proof}

As a consequence, we obtain:

\begin{corollary}
A graph $G$ is arboreal-explainable if and only if $G$ is hole-free and does not contain the house, the gem and the domino as an induced subgraph. In particular, the property of being arboreal-explainable is hereditary.
\end{corollary}

\section{Ptolemaic supergraph}\label{sec-compl}

For $(N,t)$ a labelled arboreal network, the graph $G=\mathcal C(N,t)$ only provides partial information regarding the pairs of vertices sharing an ancestor in $N$. Indeed, if $x,y \in L(N)$ are joined by an edge in $G$, then by definition, $x$ and $y$ share an ancestor in $N$. However, the converse in not necessarily true, and \new{two leaves of $N$ sharing an ancestor in $N$ may or may not be joined by an edge in $G$, depending on the label of their last common ancestor.}

Let $G=(X,E)$ be an arboreal-explainable graph, and let $(N,t)$ be a labelled arboreal network explaining $G$. Formally speaking, the above observation means that the shared ancestry graph $\mathcal A(N)$ of $N$ is a supergraph of $G$. Note that $V(\mathcal A(N))=L(N)=V(G)$ by definition. By Theorem~\ref{th:ptol}, $\mathcal A(N)$ is connected and Ptolemaic. Hence, $\mathcal A(N)$ is a \emph{Ptolemaic supergraph} (also known as \emph{Ptolemaic completion}, see \cite{CGP21}) of $G$.

\begin{figure}[h]
\begin{center}
\includegraphics[scale=0.8]{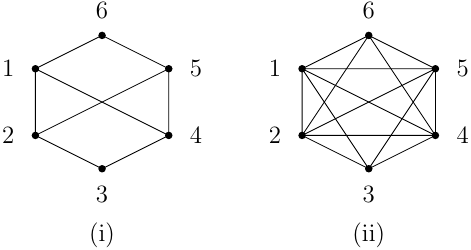}
\end{center}
\caption{(i) An undirected graph $G$ with vertex set $\{1,2,3,4,5,6\}$. (ii) A supergraph $G^*$ of $G$ on the same vertex set. Although $G$ is distance-hereditary, and therefore, arboreal-explainable, and $G^*$ is connected and Ptolemaic, there exists no labelled-arboreal network $(N,t)$ such that $\mathcal C(N,t)=G$ and $\mathcal A(N)=G^*$.}
\label{fig-com}
\end{figure}

For a pair $(G,G^*)$ of graphs such that $G$ is arboreal-explainable and $G^*$ is a connected Ptolemaic supergraph of $G$, there does not necessarily exist a labelled arboreal network $(N,t)$ such that $(N,t)$ explains $G$ and $\mathcal A(N)=G^*$. For example, consider the graph $G$ with vertex set $\{1,2,3,4,5,6\}$ depicted in Figure~\ref{fig-com}~(i). One can easily check that $G$ is distance-hereditary (and therefore, arboreal-explainable by Theorem~\ref{th:dh}), and that the graph $G^*$ depicted in Figure~\ref{fig-com}~(ii) is a Ptolemaic supergraph of $G$. However, one can verify that there is no labelled arboreal network $(N,t)$ such that $\mathcal C(N,t)=G$ and $\mathcal A(N)=G^*$.


In the following, we characterize the pairs of graphs $(G,G^*)$ for which there exists a labelled arboreal network $(N,t)$ that explains $G$ and is such that $\mathcal A(N)=G^*$. To do this, we first introduce further definitions. For $X$ and $M$ two nonempty finite sets, a \emph{symbolic map} is a map $d: {X \choose 2} \to M \cup \odot$, where $\odot \notin M$. To ease notation, for $x, y \in X$ distinct, we write $d(x,y)$ instead of $d(\{x,y\})$. In particular, $d(x,y)=d(y,x)$ always holds. Interestingly, a pair $(G,G^*)$ where $G$ is a graph and $G^*$ is a supergraph of $G$ induces a symbolic map $d_{(G,G^*)}: {V(G) \choose 2} \to \{0,1,\odot\}$, defined, for $x,y \in V(G)$ distinct, by:
\[
\begin{array}{r c  l}
d_{(G,G^*)}(x,y) = & 1 &  \text{ if } \{x,y\} \text{ is an edge of } G.\\
& 0 & \text { if } \{x,y\} \text{ is an edge of } G^* \text{ but not of } G.\\
& \odot & \text{ if } \{x,y\} \text{ is not an edge of } G^*.\\
\end{array}
\]

Following \cite{HMS24}, we say that a labelled arboreal network $(N,t)$ \emph{represents} a symbolic map $d: {X \choose 2} \to \{0,1,\odot\}$, if, for all $x,y \in X$ distinct, $x$ and $y$ share an ancestor in $N$ if and only if $d(x,y) \neq \odot$, and, if this holds, $t(\lca_N(x,y))=d(x,y)$\footnote{Note that, in \cite{HMS24} and all results therein, the map $t$ may take more than two distinct values. In particular, the definition of a labelled network given in \cite{HMS24} is more general than the one we are using here.}. As a first observation, we have:

\begin{lemma}\label{lm:diss}
Let $G=(X,E)$ be a arboreal-explainable graph, and let $G^*$ be a supergraph of $G$ with $V(G^*)=X$. Let $(N,t)$ be a labelled arboreal network with leaf set $X$. Then, the following are equivalent:
\begin{itemize}
\item[(i)] $(N,t)$ explains $G$ and $\mathcal A(N)=G^*$.
\item[(ii)] $(N,t)$ represents $d_{(G,G^*)}$.
\end{itemize}
\end{lemma}

\begin{proof}
To ease notation, we put $d=d_{(G,G^*)}$.

(i) $\Rightarrow$ (ii): Let $x,y \in X$ distinct. Suppose first that $d(x,y)=\odot$. By definition of $d$, this means that $\{x,y\}$ is not an edge of $G^*$. Since $G^*=\mathcal A(N)$, it follows that $x$ and $y$ do not share an ancestor in $N$. Suppose now that $d(x,y) \neq \odot$. In particular, $\{x,y\}$ is an edge of $G^*=\mathcal A(N)$, so $\lca_N(x,y)$ is well-defined. If $d(x,y)=1$, then by definition of $d$, $\{x,y\}$ is an edge of $G$. Since $(N,t)$ explains $G$, $t(\lca_N(x,y))=1$ follows.  If $d(x,y)=0$, then by definition of $d$, $\{x,y\}$ is not an edge of $G$. Since $(N,t)$ explains $G$, $t(\lca_N(x,y))=0$ must hold. Hence, $(N,t)$ represents $d$, so (ii) holds.

(ii) $\Rightarrow$ (i): Let $x,y \in X$ distinct. Suppose first that $\{x,y\}$ is an edge of $G$. By definition of $d$, $d(x,y)=1$. Since $(N,t)$ represents $d$, $\lca_N(x,y)$ exists, and satisfies $t(\lca_N(x,y))=1$. Suppose now that $\{x,y\}$ is not an edge of $G$. If $\{x,y\}$ is not an edge of $G^*$, then $d(x,y)=\odot$ so $x$ and $y$ do not share an ancestor in $N$. If otherwise,  $\{x,y\}$ is an edge of $G^*$, $d(x,y)=0$, so $\lca_N(x,y)$ exists, and satisfies $t(\lca_N(x,y))=0$. In summary, we have:
\begin{itemize}
\item[-] $\lca_N(x,y)$ exists and satisfies $t(\lca_N(x,y))=1$ if and only if $\{x,y\}$ is an edge of $G$.
\item[-] $\lca_N(x,y)$ exists if and only if $\{x,y\}$ is an edge of $G^*$.
\end{itemize}
The first statement implies that $(N,t)$ explains $G$, and the second statement implies that $\mathcal A(N)=G^*$. Hence, (i) holds.
\end{proof}

For $d: {X \choose 2} \to M \cup \odot$ a symbolic map, we denote by $G_d$ the undirected graph with vertex set $X$, where two distinct vertices $x,y \in X$ are joined by an edge in $G$ if $d(x,y) \neq \odot$. Symbolic maps that can be represented by a labelled arboreal network were characterized in \cite{HMS24}:

\begin{theorem}[\cite{HMS24}, Theorem~7.5]\label{th:diss}
Suppose that $d: {X \choose 2} \to M^{\odot}$ is a symbolic map. Then there exists a labelled arboreal network $(N,t)$ representing $d$ if and only if the following four properties hold:
\begin{itemize}
\item[(A1)] $G_d$ is connected and Ptolemaic.
\item[(A2)] If $x,y,z \in X$ are pairwise distinct and are such that $d(x,y), d(x,z)$ and $d(y,z)$ are pairwise distinct, then $\odot \in \{d(x,y), d(x,z),d(y,z)\}$.
\item[(A3)] If $x,y,z,u \in X$ are pairwise distinct and are such that $d(x,y)=d(y,z)=d(z,u) \neq d(y,u)=d(u,x)=d(x,z)$, then $\odot \in \{d(x,y), d(x,z)\}$. 
\item[(A4)] If $x,y,z,u \in X$ are pairwise distinct and are such that $d(z,u)= \odot$ and $d$ maps all other pairs of elements of $\{x,y,z,u\}$ to an element of $M$, then $d(x,z)=d(y,z)$ and $d(x,u)=d(y,u)$ hold.
\end{itemize}
\end{theorem}

To state the next result, we require further terminology. Let $G, G^*$ be two graphs such that $G^*$ is a supergraph of $G$. Let $(x,y,z,u)$ be an ordered set of four distinct elements of $V(G)$. We say that $(x,y,z,u)$ is an \emph{asymmetric diamond} of $(G,G^*)$ if the following three properties hold:
\begin{itemize}
\item[-] Any two elements of $\{x,y,z,u\}$ are joined by an edge in $G^*$, except $z$ and $u$.
\item[-] $x$ and $z$ are joined by an edge in $G$.
\item[-] $y$ and $z$ are not joined by an edge in $G$.
\end{itemize}

\begin{figure}[h]
\begin{center}
\includegraphics[scale=0.8]{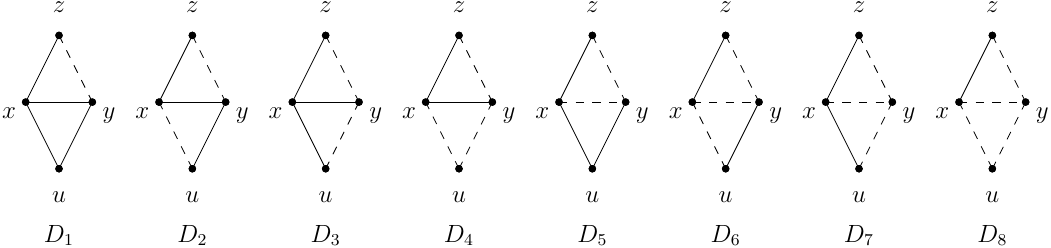}
\end{center}
\caption{The eight types of asymmetric diamonds. For $G$ a graph, $G^*$ a supergraph of $G$, and $(x,y,z,u)$ an asymmetric diamond of $(G,G^*)$, there exists exactly one $i \in\{1, \ldots, 8\}$ such that the edge set of $G[\{x,y,z,u\}]$ is precisely the set of plain edges of $D_i$. In particular, the dashed edges of $D_i$ are the edges of $G^*[\{x,y,z,u\}]$ that are not in $G$.}
\label{fig-loz}
\end{figure}

As can be seen in Figure~\ref{fig-loz}, there are $8$ distinct types of asymmetric diamonds, depending on which of the pairs $\{x,u\}, \{y,u\}$ and $\{x,y\}$ are edges of $G$.

Putting Lemma~\ref{lm:diss} and Theorem~\ref{th:diss} together, we obtain:

\begin{theorem}\label{th:ptext}
Let $G=(X,E)$ be a connected arboreal-explainable graph, and let $G^*$ be a supergraph of $G$ with $V(G^*)=X$. Then, there exists a labelled arboreal network $(N,t)$ explaining $G$ such that $\mathcal A(N)=G^*$ if and only if the following three properties hold: 
\begin{itemize}
\item[(E1)] $G^*$ is connected and Ptolemaic.
\item[(E2)] There are no four elements $x,y,z,u \in X$ such that $\{x,y,z,u\}$ induces a $P_4$ in $G$ and a clique in $G^*$.
\item[(E3)] The pair $(G,G^*)$ does not induce an asymmetric diamond.
\end{itemize}
\end{theorem}

\begin{proof}
Suppose first that there exists a labelled arboreal network $(N,t)$ explaining $G$ such that $\mathcal A(N)=G^*$. By Theorem~\ref{th:ptol}, $G^*$ is connected and Ptolemaic, so (E1) holds.

Consider now the map $d=d_{(G,G^*)}$. By Lemma~\ref{lm:diss}, $(N,t)$ represents $d$, so $d$ is a symbolic arboreal map. By Theorem~\ref{th:diss}, it follows that $d$ satisfies properties (A1) to (A4).

To see that (E2) holds, let $x,y,z,u \in X$ be such that $\{x,y,z,u\}$ induces a $P_4$ in $G$, with edges $\{x,y\}$, $\{y,z\}$ and $\{z,u\}$. We next show that $\{x,y,z,u\}$ does not induce a clique in $G^*$. By definition of $d$, we have $d(x,y)=d(y,z)=d(z,u)=1$, $d(x,z) \neq 1$, $d(x,u) \neq 1$, and $d(y,u) \neq 1$. Since $d$ satisfies (A3), $d(x,z)=d(x,u)=d(y,u)=0$ cannot hold. Hence, at least one of $d(x,z)=\odot$, $d(x,u)=\odot$ or $d(y,u)=\odot$ must hold. By definition of $d$, it follows that (at least) one of the pairs $\{x,z\}$, $\{x,u\}$ or $\{y,u\}$ is not an edge of $G^*$. In particular, $\{x,y,z,u\}$ does not induce a clique in $G^*$, so (E2) holds.

To see that (E3) holds, suppose that $x,y,z,u \in X$ are such that all pairs of elements of $\{x,y,z,u\}$ except the pair $\{z,u\}$ are edges of $G^*$, and suppose that $\{x,z\}$ is an edge of $G$. \new{By definition, $(x,y,z,u)$ is an asymetric diamond of $(G,G^*)$ if and only if $\{y,z\}$ is not an edge of $G$. Hence,} we need to show that $\{y,z\}$ is an edge of $G$. By definition of $d$, we have $d(z,u)=\odot$, while $d$ maps all other pairs of elements of $\{x,y,z,u\}$ to an element in $\{0,1\}$. Since $d$ satisfies (A4), it follows that $d(x,z)=d(y,z)$ must hold. By assumption, $\{x,z\}$ is an edge of $G$, so $d(x,z)=1$. Hence, $d(y,z)=1$, which, by definition of $d$, implies that $\{y,z\}$ is an edge of $G$, as desired.

Conversely, suppose that properties (E1), (E2), (E3) hold. By Lemma~\ref{lm:diss}, to show that there exists a labelled arboreal network $(N,t)$ explaining $G$ such that $\mathcal A(N)=G^*$, it suffices to show that $d=d_{(G,G^*)}$ is a symbolic arboreal map. Using Theorem~\ref{th:diss}, we do this by showing that $d$ satisfies properties (A1) to (A4).

To see that $d$ satisfies (A1), it suffices to remark that $G_d=G^*$. Indeed, \new{we have $V(G_d)=V(G^*)=X$ and,} by definition of $G_d$, $\{x,y\}$ is an edge of $G_d$ if and only if $d(x,y) \neq \odot$. \new{Moreover,} by definition of $d$, $d(x,y)=\odot$ if and only if $\{x,y\}$ is not an edge of $G^*$. Since $G^*$ is connected and Ptolemaic by (E1), it follows that $G_d$ is connected and Ptolemaic, so (A1) holds.

To see that $d$ satisfies (A2), it suffices to remark that the image set of $d$ is $\{0,1,\odot\}$. In particular, if $x,y,z \in X$ are such that $d(x,y)$, $d(x,z)$ and $d(y,z)$ are pairwise distinct, then exactly one of $d(x,y)=\odot$, $d(x,z)=\odot$ or $d(y,z)=\odot$ must hold.

To see that $d$ satisfies (A3), let $x,y,z,u \in X$ distinct be such that $d(x,y)=d(y,z)=d(z,u) \neq d(y,u)=d(u,x)=d(x,z)$. We \new{need to} show that one of $d(x,y)=\odot$ or $d(x,z)=\odot$ must hold. If neither $d(x,y)=1$ nor $d(x,z)=1$ holds, then since the image set of $d$ is $\{0,1,\odot\}$, one of $d(x,y)=\odot$ or $d(x,z)=\odot$ must hold. If otherwise, one of $d(x,y)=1$ or $d(x,z)=1$ holds, say $d(x,y)=1$, then by definition of $d$, $\{x,y,z,u\}$ induces a $P_4$ in $G$. By (E2), $\{x,y,z,u\}$ does not induce a clique in $G^*$, so at least one of $\{x,z\}$, $\{x,u\}$ or $\{y,u\}$ is not an edge of $G^*$. By definition of $d$, this means that at least one of $d(x,z)=\odot$, $d(x,u)=\odot$ or $d(y,u)=\odot$ holds. Since  $d(y,u)=d(u,x)=d(x,z)$ holds by choice of $\{x,y,z,u\}$, $d(x,z)=\odot$ follows.

To see that $d$ satisfies (A4), let $x,y,z,u \in X$ distinct be such that  $d(z,u)=\odot$, while $d$ maps all other pairs of elements of $\{x,y,z,u\}$ to an element in $\{0,1\}$. By definition of $d$, this means that $\{z,u\}$ is the only pair of elements of $\{x,y,z,u\}$ that is not an edge of $G^*$. \new{We need to show that $d(x,z)=d(y,z)$ and $d(x,u)=d(y,u)$ both hold. To see that $d(x,z)=d(y,z)$ holds,} suppose for contradiction that \new{this is not the case, that is,} $d(x,z) \neq d(y,z)$. Since neither $d(x,z)=\odot$ nor $d(y,z)=\odot$ holds, it follows that, up to a permutation in the set $\{x,y\}$, we have $d(x,z)=1$ and $d(y,z)=0$ .By definition of $d$, this means that $\{x,z\}$ is an edge of $G$, while $\{y,z\}$ is an edge of $G^*$  that is not an edge of $G$. In particular, $(x,y,z,u)$ is an asymmetric diamond of $(G,G^*)$. By (E3), this is impossible. Hence, we have $d(x,z)=d(y,z)$. The roles of $z$ of $u$ being symmetric, we \new{can use similar arguments to show that} $d(x,u)=d(y,u)$ also holds. Hence, $d$ satisfies (A4).
\end{proof}

\section{2-rooted arboreal networks}\label{sec-2r}

As stated in the introduction, given an unrooted graph $G$, there exists a labelled phylogenetic tree $(T,t)$ explaining $G$ if and only if $G$ is a cograph. By definition, a labelled phylogenetic tree is a labelled arboreal network with one root, so a graph $G$ can be explained by a labelled arboreal network with one root if and only if $G$ is a cograph. In this section, we go one step further and characterize graphs $G$ that can be explained by a labelled arboreal network with two roots.

To do this, we use some key properties of \gatex graphs (\cite{HS22,HS23b,HS23}), a newly introduced graph class which, as we shall see, bears similarities with arboreal-explainable graphs. We start by introducing some terminology.

We say that a phylogenetic network $N$ is a \emph{galled-tree} (sometimes also called \emph{level-1 network}, see \emph{e. g.} \cite{G03,HSS22,RV09}) if \new{all vertices of $N$ have indegree $2$ or less, and all biconnected components of $N$ contain at most one vertex of indegree $2$. As remarked in \cite{HS18}, for a biconnected component $C$ of a galled tee $N$, there exists a unique pair of vertices $r(C),\eta(C) \in V(N)$ such that $r(C)$ is an ancestor of $\eta(C)$ in $N$, there exists two disjoint paths $P^1$, $P^2$ from $r(C)$ to $\eta(C)$ in $N$ with $V(P^1) \cap V(P^2)=\{r(C),\eta(C)\}$, and $C$ is precisely the graph union of $P^1$ and $P^2$.} Following \cite{HS18}, we call a biconnected component $C$ of a galled-tree $N$ a \emph{cycle} of $N$. \new{Moreover, we call a cycle $C$} \emph{weak} if one of $P^1$ or $P^2$ has length $1$, or if both $P^1$ and $P^2$ have length $2$. We say that $C$ is \emph{well-proportioned} if one of $P^1$ or $P^2$ has length \new{$5$} or more, or if both $P^1$ and $P^2$ have length $3$ or more.

Interestingly, if $N$ is a galled-tree, then for any two leaves $x,y$ of $N$, $\lca(x,y)$ exists and is unique \cite{HSS22,HS18}. Based on this property, we say that an undirected graph $G$ with vertex set $X$ is a \emph{\gatex graph} if there exists a labelled galled-tree $(N,t)$ with leaf set $X$ such that, for all $x,y \in X$ distinct, $t(\lca_N(x,y))=1$ if and only if $\{x,y\}$ is an edge of $G$. In that case, we say that $(N,t)$ \emph{explains} $G$. We have:

\begin{lemma}[\cite{HS22}, Lemma~5.4]\label{lm:weak}
Let $G$ be an undirected graph, and let $(N,t)$ be a labelled galled-tree explaining $G$. \new{If $N$ contains a weak cycle, then} there exists a labelled galled-tree $(N',t')$ explaining $G$ such that $N'$ has strictly less cycles than $N$.
\end{lemma}

\new{For $N$ a galled tree, $C$ a cycle of $N$, and $v$ a vertex of $C$ distinct from $r(C)$ and $\eta(C)$, we have by definition that $v$ has exactly one child $v'$ in $V(C)$. In the following, we denote by $L^-_C(v)$ the set of leaves of $N$ that are descendant of $v$ but not of $v'$, that is, $L^-_C(v)=L_v \setminus L_{v'}$. The next result can be used to identify the last common ancestor of two leaves in certain situations.

\begin{lemma}\label{lm:lcaC}
Let $N$ be a galled tree, $C$ be a cycle of $N$, and $u,v$ be two distinct vertices of $V(C) \setminus \{r(C)\}$. We have:
\begin{itemize}
\item[(i)] If $u$ is an ancestor of $v$, then for all $x \in L^-_C(u), y \in L_v$, $\lca_N(x,y)=u$.
\item[(ii)] If $u$ and $v$ are incomparable, then for all $x \in L^-_C(u), y \in L^-_C(v)$, $\lca_N(x,y)=r(C)$.
\end{itemize}
\end{lemma}

\begin{proof}
(i) First, we remark that if $u$ is an ancestor of $v$, and both $u$ and $v$ are vertices of $C$, $u$ is distinct from $\eta(C)$, so $L^-_C(u)$ is well defined. So, let $x \in L^-_C(u), y \in L_v$.

Since $u$ is an ancestor of $v$, $u$ is an ancestor of both $x$ and $y$. Hence, we need to show that no child $w$ of $u$ is an ancestor of both $x$ and $y$. So, suppose for contradiction that such a vertex $w$ exists. First, we remark that by definition of $L^-_C(u)$, $w$ is not a vertex of $C$. Since both $w$ and $v$ are ancestors of $y$, there exists a path $P$ from $u$ to $y$ containing $w$, and a path $P'$ from $u$ to $y$ containing the subpath of $C$ from $u$ to $v$. Since, as remarked above, $w$ is not a vertex of $C$, these two paths are distinct. Now, consider the (necessarily unique) vertex $h \in V(P) \cap V(P')$ such that $h$ is distinct from $u$ and no proper ancestor of $h$ in $N$ other than $u$ is an element of $V(P) \cap V(P')$. It is then an easy task to verify that the graph $B$ with vertex set $V(C) \cup V(P) \cup V(P')$ and arc set $E(C) \cup E(P) \cup E(P')$ is a biconnected subgraph of $N$. But this is impossible, since $C$ is a proper subgraph of $B$ and a biconnected component of $N$. This concludes the proof that no child $w$ of $u$ is an ancestor of both $x$ and $y$.

(ii) First, we remark that if $u$ and $v$ are incomparable, and both are vertices of $C$, they are both distinct from $\eta(C)$, so $L^-_C(u)$ and $L^-_C(v)$ are well defined. So, let $x \in L^-_C(u), y \in L^-_C(v)$.

Put $r=r(C)$. Since $r$ is an ancestor of both $u$ and $v$, $r$ is an ancestor of both $x$ and $y$. Hence, we need to show that no child $w$ of $r$ is an ancestor of both $x$ and $y$. Let $u'$, $v'$ be the children of $r$ in $C$ on the paths from $r$ to $u$ and from $r$ to $v$, respectively. Since $u$ and $v$ are incomparable, $u'$ and $v'$ are distinct. The roles of $u$ and $v$ being symmetric, we may assume up to a permutation that $w$ is distinct from $v'$. Since both $w$ and $v$ are ancestors of $y$, there exists a path $P$ from $r$ to $y$ containing $w$, and a path $P'$ from $u$ to $y$ containing the subpath of $C$ from $r$ to $v$. Since $w$ is the child of $r$ in $P$, $v'$ is the child of $r$ in $P'$, and these two vertices are distinct, it follows that the paths $P$ and $P'$ are distinct. Now, consider the (necessarily unique) vertex $h \in V(P) \cap V(P')$ such that $h$ is distinct from $r$ and no proper ancestor of $h$ in $N$ other than $r$ is an element of $V(P) \cap V(P')$. It is then an easy task to verify that the graph $B$ with vertex set $V(C) \cup V(P) \cup V(P')$ and arc set $E(C) \cup E(P) \cup E(P')$ is a biconnected subgraph of $N$. But this is impossible, since $C$ is a proper subgraph of $B$ and a biconnected component of $N$. This concludes the proof that no child $w$ of $r$ is an ancestor of both $x$ and $y$.
\end{proof}}

\new{Next,} we say that a galled-tree $N$ with root $\rho$ is a \emph{basic galled-tree} if $N$ has exactly one cycle $C$, $r(C)=\rho$, and $\rho$ has outdegree $2$. Note that the latter condition implies that $\rho$ has no children outside of $C$. If $N$ is a basic galled-tree and $t$ is a labelling of $N$ satisfying $t(\rho)=0$ (\emph{resp.} $t(\rho)=1$), we say that $(N,t)$ is a \emph{$0$-basic labelled galled-tree} (\emph{resp.} a \emph{$1$-basic labelled galled-tree})

As a first observation, we have:

\begin{lemma}\label{lm:2rgt}
Let $G$ be an undirected graph. There exists a labelled arboreal network $(N,t)$ with exactly two roots explaining $G$ if and only if there exists a $0$-basic labelled galled-tree $(N',t')$ explaining $G$.
\end{lemma}

\begin{proof}
Suppose first that there exists a labelled arboreal network $(N,t)$ with exactly two roots explaining $G$. Let $\rho_1$ and $\rho_2$ be the roots of $N$. Consider the phylogenetic network $N'$ obtained from $N$ by first adding a new vertex $\rho$, and adding arcs $(\rho,\rho_1)$ and $(\rho,\rho_2)$. \new{Note that $V^*(N')=V^*(N) \cup \{\rho\}$.} For all $v \in V^*(N')$, define $t'(v)=t(v)$ if $v \in V^*(N)$, and $t'(v)=0$ if $v=\rho$.

\new{We first show that $(N',t')$ is a $0$-basic galled-tree. By Lemma~\ref{lm:rh}, we have $\sum\limits_{h \in H(N)} \mathrm{indeg}(h)-1=|R(N)|-1=1$, so $N$ has exactly one vertex $h$ of indegree $2$ or more, and $h$ has indegree exactly $2$. Moreover, since $N$ is arboreal, $h$ is a descendant of both $\rho_1$ and $\rho_2$. Consider now the subset $V'$ of $V(N')$ containing all vertices of $N$ that are ancestors of $h$, and let $C$ be the subgraph of $N'$ induced by the elements of $N'$. Note that $C$ consist of the path from $\rho_1$ to $h$, the path from $\rho_2$ to $h$, and the arcs $(\rho,\rho_1)$ and $(\rho,\rho_2)$. In particular, $C$ is biconnected. Moreover, by definition of $V'$ and uniqueness of $h$, for all arcs $(u,v)$ of $N'$ not in $C$, all descendants of $v$ have indegree $1$ in $N'$. In particular, the removal of $v$ disconnectes $N'$, so $(u,v)$ does not belong to a biconnected component of $N'$. Hence, $C$ is a biconnected component of $N'$, and $N'$ does not have any biconnected component other than $C$. Since $C$ contains exactly one vertex of indegree $2$ (the vertex $h$), $N'$ is a galled-tree. Moreover, we have $r(C)=\rho$, and $t'(\rho)=0$, so} $(N',t')$ is a $0$-basic labelled galled-tree It remains to show that $(N',t')$ explains $G$.

\new{To see this, we remark that for all vertices $v$ of $N$, the set of leaves that are descendants of $v$ in $N$ coincides with the set of leaves that are descendants of $v$ in $N'$. In particular, $\lca_{N'}(x,y)=\lca_{N}(x,y)$ holds for all $x,y \in V(G)$ such that $x$ and $y$ share an ancestor in $N$. If otherwise, $x$ and $y$ are such that $x$ and $y$ do not share an ancestor in $N$, then $\lca_{N'}(x,y)=\rho$. Since $t'(v)=t(v)$ for all $v \in V^*(N)$, and $t'(\rho)=0$, and since $(N,t)$ explains $G$, it follows that} $\{x,y\}$ is an edge of $G$ if and only if $t'(\lca_{N'}(x,y))=1$. This concludes the proof that $(N',t')$ explains $G$.

Conversely, suppose that there exists a $0$-basic labelled galled-tree $(N',t')$ explaining $G$. Let $\rho$ be the root of $N'$, and consider the graph $N$ obtained from $N'$ by removing $\rho$ and its incident arcs. Since $N'$ is basic, $N$ has a unique cycle $C$, and we have $\rho=r(C)$. In particular, $\rho$ belongs to a biconnected component of $N'$, so the graph $N$ is connected. Thus, $N$ is a network.
\new{Moreover, since $\rho$ has outdegree $2$ in $N'$ by assumption, $N$ has exactly two roots, namely, the children of $\rho$ in $N'$. Finally, since $N$ is a basic galled-tree, $N'$ has exactly one vertex $h$ of indegree $2$ or more, and that vertex has indegree exactly $2$. Note that except for $\rho_1$ and $\rho_2$, all vertices of $N$ have the same indegree in $N'$, so $h$ is also the only vertex of $N$ of indegree $2$ or more, and has indegree exactly $2$. Hence, we have $\sum\limits_{h \in H(N)} \mathrm{indeg}(h)-1=1=|R(N)|-1$, from which Lemma~\ref{lm:rh} implies that $N$ is arboreal. In summary,} $N$ is an arboreal network with exactly two roots. It remains to show that $(N,t)$ explains $G$, where $t$ is the restriction of $t'$ to $V^*(N)=V^*(N') \setminus \{\rho\}$.

\new{To see this, we remark that for all vertices $v$ of $N$, the set of leaves that are descendants of $v$ in $N$ coincides with the set of leaves that are descendants of $v$ in $N'$. In particular, $\lca_{N'}(x,y)=\lca_{N}(x,y)$ holds for all $x,y \in V(G)$ such that $x$ and $y$ share an ancestor in $N$. If otherwise, $x$ and $y$ are such that $x$ and $y$ do not share an ancestor in $N$, then $\lca_{N'}(x,y)=\rho$. Since $t'(v)=t(v)$ for all $v \in V^*(N)$, and $t'(\rho)=0$, and since $(N',t')$ explains $G$, it follows that $\{x,y\}$ is an edge of $G$ if and only if $x$ and $y$ share an ancestor in $N$, and $t(\lca_{N}(x,y))=1$.} This concludes the proof that $(N,t)$ explains $G$.
\end{proof}

\begin{figure}[h]
\begin{center}
\includegraphics[scale=1]{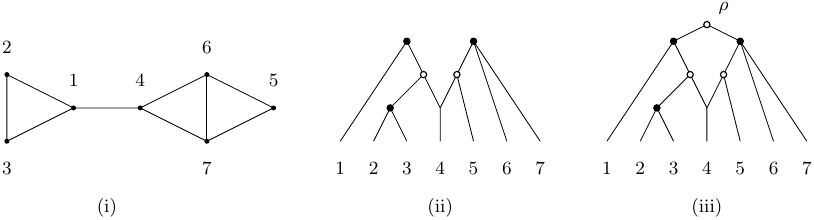}
\end{center}
\caption{(i) An undirected graph $G$ with vertex set $\{1,2,3,4,5,6,7\}$. (ii) A labelled arboreal network with two roots explaining $G$. (iii) A $0$-basic labelled galled-tree explaining $G$. In both (ii) and (iii), vertices with label $1$ are indicated with $\bullet$, and vertices with label $0$ are indicated with $\circ$.}
\label{fig-2rgt}
\end{figure}

To illustrate Lemma~\ref{lm:2rgt}, consider the graph $G$ depicted in Figure~\ref{fig-2rgt}~(i). A labelled arboreal network $(N,t)$ with two roots explaining $G$ is presented in Figure~\ref{fig-2rgt}~(ii), and a $0$-basic labelled galled-tree $(N',t')$ explaining $G$ is presented in Figure~\ref{fig-2rgt}~(iii). One can see that $(N',t')$ is obtained from $(N,t)$ by adding a new vertex $\rho$ as the parent of both roots of $N$, and defining $t(\rho)=0$, which is precisely the construction described in the first part of the proof of Lemma~\ref{lm:2rgt}. Conversely, $(N,t)$ is obtained from $(N',t')$ by removing the vertex $\rho$, which is the operation used in the second part of the proof of Lemma~\ref{lm:2rgt}.

In the following, we focus on graphs $G$ that can be explained by a labelled arboreal network with exactly two roots, but cannot be explained by a labelled arboreal network with \new{one} root. In particular, we ignore the special case where $G$ is a cograph. We have:

\begin{proposition}\label{pr:2rgt}
Let $G$ be an undirected graph that is not a cograph. There exists a labelled arboreal network $(N,t)$ with exactly two roots explaining $G$ if and only if $G$ is a distance-hereditary \gatex graph, and exactly one connected component of $G$ is not a cograph.
\end{proposition}

\begin{proof}
Suppose first that there exists a labelled arboreal network $(N,t)$ with exactly two roots explaining $G$. By Theorem~\ref{th:dh}, $G$ is distance-hereditary. Moreover, by Lemma~\ref{lm:2rgt}, there exists a labelled galled-tree $(N',t')$ explaining $G$, so $G$ is a \gatex graph.

It remains to show that exactly one connected component of $G$ is not a cograph. Since $G$ is not a cograph, and the disjoint union of cographs is a cograph, at least one connected component of $G$ is not a cograph. To see that there cannot be more than one such component, suppose for contradiction that $G$ has two distinct connected components $G_1$ and $G_2$ such that neither $G_1$ nor $G_2$ are cographs. Note that since $N$ is arboreal and has two roots, it follows from Lemma~\ref{lm:rh} that $H(N)$ contains exactly one element, which we call $h$. By Lemma~\ref{lm:indcog}, there exists $x_1, y_1 \in V(G_1)$ distinct such that $x_1 \in L_h$, $y_1 \notin L_h$, and $\{x_1,y_1\}$ is an edge of $G$. Similarly, there exists $x_2, y_2 \in V(G_2)$ distinct such that $x_2 \in L_h$, $y_2 \notin L_h$, and $\{x_2,y_2\}$ is an edge of $G$. Note that by choice of $G_1$ and $G_2$, $x_1$ and $x_2$ are distinct. \new{Since $(N,t)$ explains $G$, $x_1$ and $y_1$ share an ancestor in $N$, and we have $t(v)=1$ for $v=\lca_N(x_1,y_1)$. Next, we remark that since $x_1 \in L_h$, $y_1 \notin L_h$, and $h$ is the only vertex of $N$ of indegree $2$ or more, $v$ must be a strict ancestor of $h$. In particular, $v$ is also an ancestor of $x_2$.

We next show that $v=\lca_N(x_2,y_1)$. To see this, suppose for contradiction that $v \neq \lca_N(x_2,y_1)$. Then, there exists a child $v'$ of $v$ in $N$ such that both $x_2$ and $y_1$ are descendant of $v'$. Moreover, $v=\lca_N(x_1,y_1)$, so $v'$ is not an ancestor of $x_1$. Hence, there exists a path $P$ from $v$ to $x_2$ that does not contain $h$. Since $x_2 \in L_h$, and $v$ is an ancestor of $h$ there also exists a path $P'$ from $v$ to $x_2$ that contains $h$. Now, let $w$ be an element of $V(P) \cap V(P')$ such that no parent of $w$ is an element of $V(P) \cap V(P')$. By choice of $P$ and $P'$, $w$ is a strict descendant of $h$. In particular, $w$ is distinct from $v$, so $w$ has one parent in $P$ and one parent in $P'$, and these parents are distinct. Hence, $w \in H(N)$, a contradiction since $w \neq h$ and $h$ is the only element of $H(N)$. Therefore, $v=\lca_N(x_2,y_1)$ as desired. Since $t(v)=1$ and $(N,t)$ explains $G$, it follows} that $\{x_2,y_1\}$ is an edge of $G$. This is a contradiction, since $x_2 \in V(G_2)$, $y_1 \in V(G_1)$, and $G_1$ and $G_2$ are distinct connected components of $G$.

Conversely, suppose that $G$ is distance-hereditary and \gatex, and that exactly one connected component $G'$ of $G$ is not a cograph. We first show that there exists a labelled arboreal network $(N^-,t^-)$ with exactly two roots explaining $G'$. In view of Lemma~\ref{lm:2rgt}, it suffices to show that there exists a $0$-basic labelled galled-tree $(N',t')$ explaining $G'$.

Since $G$ is \gatex, and the property of being \gatex is hereditary, $G'$ is \gatex. Hence, there exists a labelled galled-tree $(N',t')$ explaining $G'$. Without loss of generality, we may assume that the number of cycles in $N'$ is minimum. \new{In particular, Lemma~\ref{lm:weak} implies that $N'$ does not contain any weak cycle.} Note that since $G'$ is not a cograph, $N'$ is not a phylogenetic tree, so $N'$ has at least one cycle. We may also assume that (i) there are no cycle $C$ and no $u,v \in V^*(N') \cap V(C)$ distinct such that $(u,v)$ is an arc of $N'$, and $t'(u)=t'(v)$ (see \cite[Proposition~5.3]{HS22}, which states that if one such pair of vertices exists, then identifying the vertices $u$ and $v$ results in a labelled galled-tree still explaining $G'$) and (ii) that all vertices of a cycle $C$ have outdegree $2$ \new{in $N$}, except $\eta(C)$ which has outdegree $1$ (see \cite[Definition 1 and Observation 1]{HS23b}, which shows that it is always possible to build a labelled galled-tree explaining $G'$ that enjoys this property).

We first show that $N'$ is a basic galled-tree. Since $G'$ is not a cograph, $N'$ is not a tree, so $N'$ contains at least one cycle. Let $C$ be one such cycle, and let $v=r(C)$. If $G'[L_v]$ is a cograph, then one can replace the labelled subnetwork of $N'$ rooted at $v$ with a labelled phylogenetic tree explaining $G[L_v]$. Clearly, the resulting network is a galled-tree explaining $G$, and that galled-tree has strictly less cycles than $N'$. This is impossible, since the number of cycles in $N'$ is minimum by assumption. Hence, $G'[L_v]$ is not a cograph. In particular, there exists $Y \subseteq G'[L_v]$ such that $G'[Y]$ is a $P_4$. Suppose now for contradiction that $v$ is distinct from $\rho$. In particular, $G'[L_v]$ is a proper subgraph of $G'$. If all proper ancestors of $v$ in $N'$ have label $0$, then $G'[L_v]$ is a connected component of $G'$ distinct from $G'$, which is impossible since $G'$ is connected. Hence, $v$ has an ancestor $u$ in $N'$ such that \new{$u \in V^*(N')$ and} $t'(u)=1$. In this case, consider a child $u'$ of $u$ that is not an ancestor of $v$, and let $x \in L_{u'}$. By choice of $u'$, we have that $\lca_{N'}(x,y)=u$ for all $y \in L_v$. Since $t(u)=1$ and $(N',t')$ explains $G'$, it follows that $\{x,y\}$ is an edge of $G'$ for all $y \in L_v$. In particular, $\{x,y\}$ is an edge of $G'$ for all $y \in Y$, so $G'[Y \cup \{x\}]$ is a gem. This is impossible, since $G$ is distance-hereditary, and therefore, does not have a gem as an induced subgraph. Hence, all cycles $C$ of $N$ satisfy $r(C)=\rho$. By choice of $N'$, $\rho=r(C)$ has outdegree $2$. In particular, since all cycles $C$ of $N'$ satisfy $r(C)=\rho$, this implies that $N'$ contains exactly one cycle, which concludes the proof that $N'$ is basic.

It remains to show that $(N',t')$ can be chosen to be $0$-basic, that is, that up to modifying $N'$ and $t'$, $t'(\rho)=0$ holds. To see this, suppose that $t'(\rho)=1$. Let $C$ be the unique cycle of $N'$, and let $\eta=\eta(C)$. \new{Recall that} by choice of $N'$, $C$ is not weak. \new{Now, let $P^1$ and $P^2$ be the two paths between $\rho$ and $\eta$ in $C$. Without loss of generality, we may assume that $|V(P^1)| \leq |V(P^2)|$.} We distinguish between two cases: (a) $C$ is well proportioned and (b) $C$ is not well-proportioned.

\begin{itemize}
\item[Case (a):] $C$ is well proportioned. \new{We show that this case is impossible. Since $C$ is not weak, we either have (a1) $|V(P^1)| \geq 4$ or (a2) $|V(P^1)|=3$ and $|V(P^2)| \geq 6$. We investigate these two subcases in turn.

\begin{itemize}
\item[(a1)] $|V(P^1)| \geq 4$. Note that by choice of $P^1$, $|V(P^2)| \geq 4$ also holds. For $i \in \{1,2\}$, let $u_i$ be the child of $\rho$ in $P^i$, and let $v_i$ be the child of $u_i$ in $P^i$. Note that since $|V(P^1)| \geq 4$ and $|V(P^2)| \geq 4$ both hold, $u_1,u_2,v_1$ and $v_2$ are well-defined, and are distinct from $\eta$. Next, for $i \in \{1,2\}$, choose $x_i \in L^-_C(u_i)$, $y_i=L^-_C(v_i)$, and $z \in L_{\eta}$. Since $u_1,u_2,v_1,v_2$ and $\eta$ are pairwise distinct, the elements $x_1,x_2,y_1,y_2$ and $z$ of $V(G')$ are pairwise distinct. By Lemma~\ref{lm:lcaC}, we have $\lca_{N'}(x_i,y_i)=\lca_{N'}(x_i,z)=u_i$, $\lca_{N'}(y_i,z)=v_i$, and $\lca_{N'}(x_1,x_2)=\lca_{N'}(x_1,y_2)=\lca_{N'}(y_1,x_2)=\lca_{N'}(y_1,y_2)=\rho$. By choice of $(N',t')$, we have $t'(\rho)=t'(v_1)=t'(v_2)=1$, and $t'(u_1)=t'(u_2)=0$. Since $(N',t')$ explains $G'$, it follows that $G'[\{x_1,x_2,y_1,y_2,z\}]$ has edges $\{x_1,x_2\}$, $\{x_1,y_2\}$, $\{y_1,x_2\}$, $\{y_1,y_2\}$, $\{y_1,z\}$, and $\{y_2,z\}$. Hence, $\{x_1,x_2,y_1,y_2,z\}$ induces a house in $G$, which is impossible since $G$ is distance-hereditary, and thus, does not have a house as an induced subgraph.

\item[(a2)]$|V(P^1)|=3$ and $|V(P^2)| \geq 6$. In this case, $P^1$ contains an unique vertex $u$ distinct from $\rho$ and $\eta$. Now, let $v_1,v_2,v_3,v_4$ be the vertices of $P^2$ such that $v_1$ is the child of $\rho$ in $P^2$, and for $j \in \{2,3,4\}$, $v_j$ is the child of $v_{j-1}$ in $P^2$. Note that since $|V(P^2)| \geq 6$, these vertices are well-defined, and are distinct from $\eta$. By choice of $N'$, we have $t'(\rho)=t'(v_2)=t'(v_4)=1$ and $t'(u)=t'(v_1)=t'(v_3)=0$. Now, let $x \in L^-_C(u)$, $z \in L_{\eta}$, and for $j \in \{1,2,3,4\}$, $y_j \in L^-_C(v_i)$. Since $u,v_1,v_2,v_3,v_4$ and $\eta$ are pairwise distinct, the elements $x,y_1,y_2,y_3,y_4$ and $z$ are also pariwise distinct. By Lemma~\ref{lm:lcaC}, we have $\lca_{N'}(x,y_j)=\rho$ and $\lca_{N'}(y_j,z)=v_j$ for all $j \in\{1,2,3,4\}$, and also $\lca_{N'}(x,z)=u$, $\lca_{N'}(y_2,y_3)=\lca_{N'}(y_2,y_4)=v_2$, and $\lca_{N'}(y_3,y_4)=v_3$. Since $(N',t')$ explains $G$, one easily verifies that these equalities imply that $\{x,y_2,y_3,y_4,z\}$ induces a gem in $G$, which is impossible since $G$ is distance-hereditary, and thus, does not have a house as an induced subgraph.}
\end{itemize}

\item[Case (b):] $C$ is not well-proportioned. \new{We show that, in that case, we can use the structure of $(N',t')$ to build a $0$-basic galled-tree $(N'',t'')$ explaining $G$. Since $C$ is not weak, we have $|V(P^1)|=3$, and $|V(P^2)| \in \{4,5\}$. We therefore} distinguish between two subcases: (b1) $|V(P^2)|=4$ and (b2) $|V(P^2)|=5$.

\begin{itemize}
\item[(b1)] $|V(P^2)|=4$. In this case, let $\{\rho,\eta,u_1,u_2,v_2\}$ be the vertex set of $C$, where $u_1, u_2$ and $v_2$ are chosen in such a way that the arc set of $C$ is $\{(\rho,u_1), (u_1,\eta), (\rho,u_2), (u_2,v_2),$ $(v_2,\eta)\}$. By choice of $(N',t')$, we have $t'(\rho)=t'(v_2)=1$, and $t'(u_1)=t'(u_2)=0$. Now, consider the network $N''$ obtained from $N'$ by removing the arcs $(u_1,\eta), (\rho,u_2)$ and $(v_2,\eta)$, and adding the arcs $(\rho,\eta), (\eta,v_2)$ and $(u_1,u_2)$ (Figure~\ref{fig-pr73}~(i)). Clearly, $N''$ is a basic galled-tree, whose unique cycle $C_0$ satisfies $r(C_0)=\rho$ and $\eta(C_0)=v_2$. Note that, since $v_2$ has outdegree $2$ in $N'$ by assumption on $N'$, $v_2$ has outdegree $1$ in $N''$. In particular, $v_2 \notin V^*(N'')$.  Now, let $t''$ be a labelling of $V^*(N'')$ defined by $t''(v)=t'(v)$ if $v \notin \{\rho, u_1, \eta\}$, $t''(\rho)=0$, and $t''(u_1)=t''(\eta)=1$. Since $t''(\rho)=0$, $(N'',t'')$ is a $0$-basic labelled galled-tree. It remains to show that $(N'',t'')$ explains $G'$. Since $(N',t')$ explains $G$, it suffices to show that $t''(\lca_{N''}(x,y))=t'(\lca_{N'}(x,y))$ holds for all $x,y \in V(G')$ distinct.

Let $x,y \in V(G'')$ distinct. If $\lca_{N'}(x,y) \notin V(C)$, then $\lca_{N''}(x,y)=\lca_{N'}(x,y)$, and\\ $t''(\lca_{N''}(x,y))=t'(\lca_{N'}(x,y))$ holds by definition of $t''$. Suppose now that $\lca_{N'}(x,y) \in V(C)$. Note that since $\eta$ has outdegree $1$ in $N'$ by assumption on $N'$, $\lca_{N'}(x,y)$ is distinct from $\eta$. If $\lca_{N'}(x,y)=u_1$, then since $u_1$ has outdegree $2$ in $N'$, we have (up to permutation) $x \in L^-_C(u_1)$ and $y \in L_{\eta}$. In that case, we have $\lca_{N''}(x,y)=\rho$ \new{by Lemma~\ref{lm:lcaC}}. Since $t''(\rho)=0=t'(u_1)$, $t''(\lca_{N''}(x,y))=t'(\lca_{N'}(x,y))$ follows. If $\lca_{N'}(x,y)=v_2$, then since $v_2$ has outdegree $2$ in $N'$, we have (up to permutation) $x \in L^-_C(v_2)$ and $y \in L_{\eta}$. In that case, we have $\lca_{N''}(x,y)=\eta$ \new{by Lemma~\ref{lm:lcaC}}. Since $t''(\eta)=1=t'(v_2)$, $t''(\lca_{N''}(x,y))=t'(\lca_{N'}(x,y))$ follows. If $\lca_{N'}(x,y)=u_2$, then since $u_2$ has outdegree $2$ in $N'$, we have (up to permutation) $x \in L^-_C(u_2)$ and $y \in L_{v_2}=L^-_C(v_2) \cup L_{\eta}$. In that case, \new{it follows from Lemma~\ref{lm:lcaC}} that we have either $\lca_{N''}(x,y)=u_2$ (if $y \in L^-_C(v_2)$), or $\lca_{N''}(x,y)=\rho$ (if $y \in L_{\eta}$). Since $t''(u_2)=t''(\rho)=0=t'(u_2)$, $t''(\lca_{N''}(x,y))=t'(\lca_{N'}(x,y))$ follows. Finally, if $\lca_{N'}(x,y)=\rho$, then we have (up to permutation) $x \in L^-_C(u_1)$ and $y \in L^-_C(u_2) \cup L^-_C(v_2)$. In that case, we have $\lca_{N''}(x,y)=u_1$ \new{by Lemma~\ref{lm:lcaC}}. Since $t''(u_1)=1=t'(\rho)$, $t''(\lca_{N''}(x,y))=t'(\lca_{N'}(x,y))$ follows.

\item[(b2)] $|V(P^2)|=5$. In this case, let $\{\rho,\eta,u_1,u_2,v_2,w_2\}$ be the vertex set of $C$, where $u_1, u_2, v_2$ and $w_2$ are chosen in such a way that the arc set of $C$ is $\{(\rho,u_1), (u_1,\eta),$ $(\rho,u_2), (u_2,v_2), (v_2, w_2), (w_2, \eta)\}$. By choice of $(N',t')$, we have $t'(\rho)=t'(v_2)=1$, and $t'(u_1)=t'(u_2)=t'(w_2)=0$. Now, consider the network $N''$ obtained from $N'$ by removing the arcs $(u_1,\eta), (\rho,u_2),$ $(u_2,v_2), (v_2, w_2)$ and $(w_2,\eta)$, and adding the arcs $(\rho,\eta), (\eta,v_2), (u_1,u_2), (u_2,w_2)$ and $(w_2,v_2)$ (Figure~\ref{fig-pr73}~(ii)). Clearly, $N''$ is a basic galled-tree, whose unique cycle $C_0$ satisfies $r(C_0)=\rho$ and $\eta(C_0)=v_2$. Note that, since $v_2$ has outdegree $2$ in $N'$ by assumption on $N'$, $v_2$ has outdegree $1$ in $N''$. In particular, $v_2 \notin V^*(N'')$.  Now, let $t''$ be a labelling of $V^*(N'')$ defined by $t''(v)=t'(v)$ if $v \notin \{\rho, u_1, w_2, \eta\}$, $t''(\rho)=0$, and $t''(u_1)=t''(w_2)=t''(\eta)=1$. Since $t''(\rho)=0$, $(N'',t'')$ is a $0$-basic labelled galled-tree. It remains to show that $(N'',t'')$ explains $G'$. Since $(N',t')$ explains $G$, it suffices to show that $t''(\lca_{N''}(x,y))=t'(\lca_{N'}(x,y))$ holds for all $x,y \in V(G')$ distinct.

Let $x,y \in V(G'')$ distinct. If $\lca_{N'}(x,y) \notin V(C)$, then $\lca_{N''}(x,y)=\lca_{N'}(x,y)$, and\\ $t''(\lca_{N''}(x,y))=t'(\lca_{N'}(x,y))$ holds by definition of $t''$. Suppose now that $\lca_{N'}(x,y) \in V(C)$. Note that since $\eta$ has outdegree $1$ in $N'$ by assumption on $N'$, $\lca_{N'}(x,y)$ is distinct from $\eta$. If $\lca_{N'}(x,y)=u_1$, then since $u_1$ has outdegree $2$ in $N'$, we have (up to permutation) $x \in L^-_C(u_1)$ and $y \in L_{\eta}$. In that case, we have $\lca_{N''}(x,y)=\rho$ \new{by Lemma~\ref{lm:lcaC}}. Since $t''(\rho)=0=t'(u_1)$, $t''(\lca_{N''}(x,y))=t'(\lca_{N'}(x,y))$ follows. If $\lca_{N'}(x,y)=w_2$, then since $w_2$ has outdegree $2$ in $N'$, we have (up to permutation) $x \in L^-_C(w_2)$ and $y \in L_{\eta}$. In that case, we have $\lca_{N''}(x,y)=\rho$ \new{by Lemma~\ref{lm:lcaC}}. Since $t''(\rho)=0=t'(w_2)$, $t''(\lca_{N''}(x,y))=t'(\lca_{N'}(x,y))$ follows. If $\lca_{N'}(x,y)=v_2$, then since $v_2$ has outdegree $2$ in $N'$, we have (up to permutation) $x \in L^-_C(v_2)$ and $y \in L_{w_2}=L^-_C(w_2) \cup L_{\eta}$. In that case, \new{it follows from Lemma~\ref{lm:lcaC}} that we have either $\lca_{N''}(x,y)=w_2$ (if $y \in L^-_C(w_2)$), or $\lca_{N''}(x,y)=\eta$ (if $y \in L_{\eta}$). Since $t''(w_2)=t''(\eta)=1=t'(v_2)$, $t''(\lca_{N''}(x,y))=t'(\lca_{N'}(x,y))$ follows. If $\lca_{N'}(x,y)=u_2$, then since $u_2$ has outdegree $2$ in $N'$, we have (up to permutation) $x \in L^-_C(u_2)$ and $y \in L_{v_2}=L^-_C(v_2) \cup L^-_C(w_2) \cup L_{\eta}$. In that case, \new{it follows from Lemma~\ref{lm:lcaC} that} we have either $\lca_{N''}(x,y)=u_2$ (if $y \in L^-_C(v_2) \cup L^-_C(w_2)$), or $\lca_{N''}(x,y)=\rho$ (if $y \in L_{\eta}$). Since $t''(u_2)=t''(\rho)=0=t'(u_2)$, $t''(\lca_{N''}(x,y))=t'(\lca_{N'}(x,y))$ follows. Finally, if $\lca_{N'}(x,y)=\rho$, then since $v_2$ has outdegree $2$ in $N'$, we have (up to permutation) $x \in L^-_C(u_1)$ and $y \in L^-_C(u_2) \cup L^-_C(v_2) \cup L^-_C(w_2)$. In that case, we have $\lca_{N''}(x,y)=u_1$ \new{by Lemma~\ref{lm:lcaC}}. Since $t''(u_1)=1=t'(\rho)$, $t''(\lca_{N''}(x,y))=t'(\lca_{N'}(x,y))$ follows.
\end{itemize}
\end{itemize}

This concludes the proof that there exists a $0$-basic labelled galled-tree $(N',t')$ explaining $G'$. By Lemma~\ref{lm:2rgt}, there exists a labelled arboreal network $(N^-,t^-)$ with exactly two roots explaining $G'$. If $G=G'$, the conclusion follows. Otherwise, since $G'$ is a connected component of $G$, there exists an induced subgraph $G''$ of $G$ such that $G$ is the disjoint union of $G'$ and $G''$. By assumption on $G$, all connected components of $G''$ are cographs, so $G''$ is a cograph. In particular, there exists a labelled phylogenetic tree $(T, t_T)$ (in other words, a labelled arboreal network with one root) explaining $G''$. By Lemma~\ref{lm:ccc}, it follows that there exists a labelled arboreal network $(N,t)$ with exactly two roots explaining $G$, which concludes the proof.
\end{proof}

\begin{figure}[h]
\begin{center}
\includegraphics[scale=0.8]{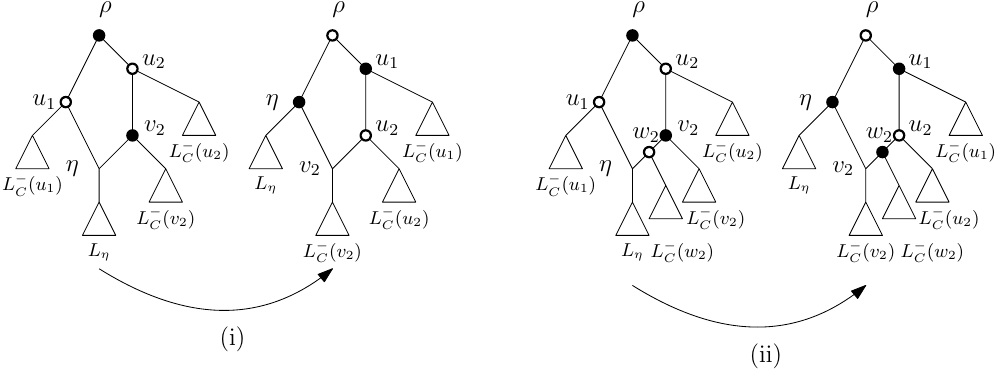}
\end{center}
\caption{(i) The transformation from $(N',t')$ to $(N'',t'')$ described in case (b1) of the proof of Proposition~\ref{pr:2rgt}. (ii) The transformation from $(N',t')$ to $(N'',t'')$ described in case (b2) of the proof of Proposition~\ref{pr:2rgt}.}
\label{fig-pr73}
\end{figure}

\begin{figure}[h]
\begin{center}
\includegraphics[scale=0.6]{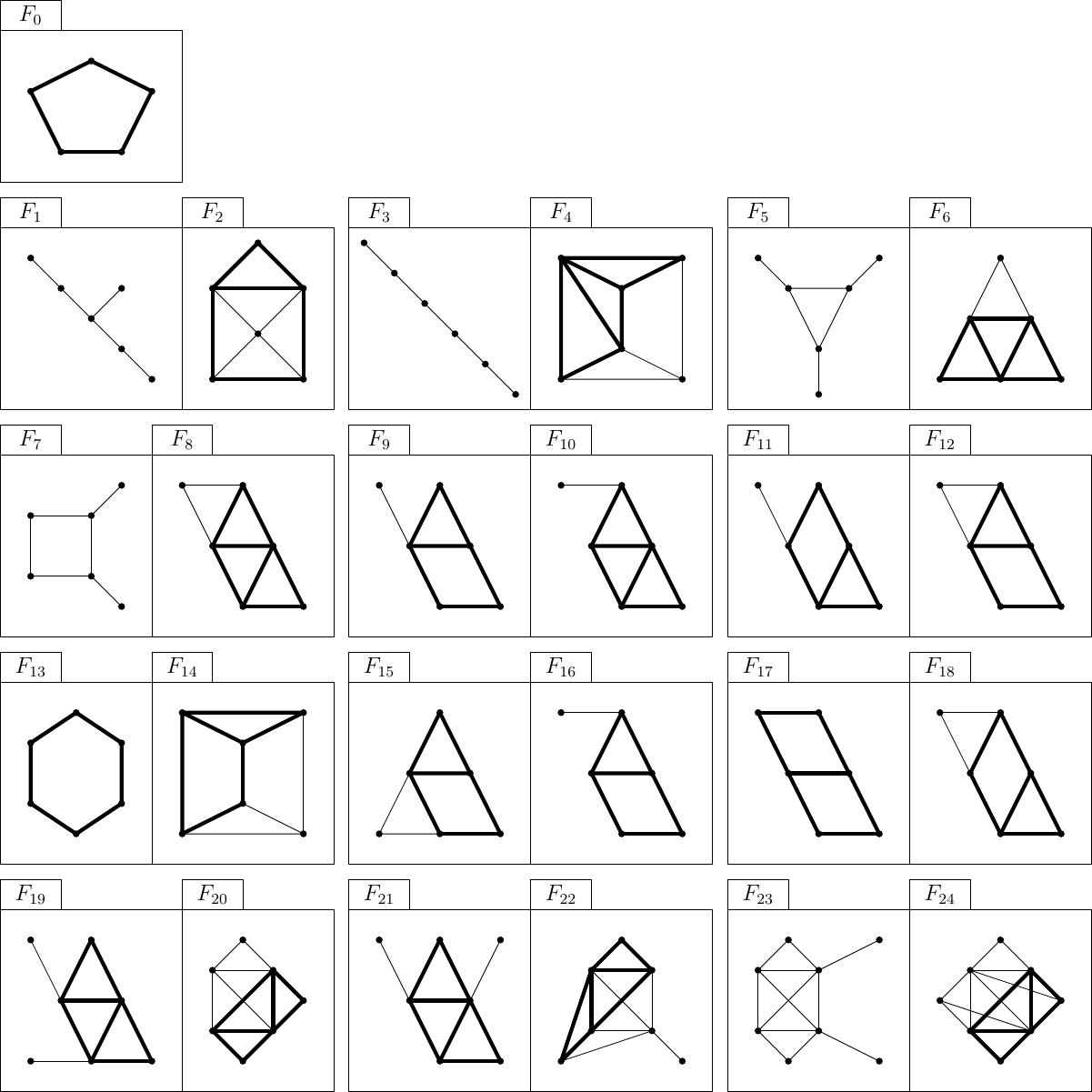}
\end{center}
\caption{The 25 forbidden induced subgraph characterizing \gatex graphs (see \cite{HS23}). In all non distance-hereditary graphs, an induced hole, house, gem, or domino is indicated in terms of bold edges.}
\label{fig-gf}
\end{figure}

In \cite{HS23}, \gatex graphs were characterized in terms of 25 forbidden induced subgraphs. These graphs are depicted in Figure~\ref{fig-gf}. Out of these 25 graphs, 5 are distance-hereditary: $F_1$, $F_3$, $F_5$, $F_7$, and $F_{23}$. All other graphs contain at least one induced hole, house, gem or domino, highlighted with bold edges in Figure~\ref{fig-gf}, and are therefore not distance-hereditary.

Putting these observations together, we obtain:

\begin{theorem}\label{th:2rgt}
Let $G$ be an undirected graph that is not a cograph. The following conditions are equivalent:
\item[(P1)] There exists a labelled arboreal network $(N,t)$ with exactly two roots explaining $G$.
\item[(P2)] There exists a $0$-basic labelled galled-tree $(N',t')$ explaining $G$.
\item[(P3)] $G$ is distance-hereditary and \gatex, and exactly one connected component of $G$ is not a cograph.
\item[(P4)] $G$ is hole-free, does not contain the house, the gem, the domino, or one of $F_1, F_3, F_5, F_7, F_{23}$ as an induced subgraph, and all induced $P_4$s of $G$ belong to the same connected component of $G$.
\end{theorem}

\begin{proof}
(P1) $\Leftrightarrow$ (P2) is stated in Lemma~\ref{lm:2rgt}, and (P1) $\Leftrightarrow$ (P3) is stated in Proposition~\ref{pr:2rgt}. We next show that (P3) $\Leftrightarrow$ (P4).

(P3) $\Rightarrow$ (P4). Since $G$ is distance-hereditary, $G$ is hole-free and does not contain the house, the gem or the domino as an induced subgraph. Moreover, $G$ is \gatex, so $G$ does not contain any of $F_1, F_3, F_5, F_7$ and $ F_{23}$ as an induced subgraph. Finally, since exactly one connected component $G'$ of $G$ is not a cograph, it follows that no connected component of $G$ distinct from $G'$ contains an induced $P_4$.

(P4) $\Rightarrow$ (P3). Since $G$ is hole-free and does not contain the house, the gem or the domino as an induced subgraph, $G$ is distance-hereditary. Moreover, the property of being distance-hereditary is hereditary, $G$ does not contain any of the non distance-hereditary forbidden subgraphs characterizing \gatex graph. Since $F_1, F_3, F_5, F_7$ and $F_{23}$ are the only distance-hereditary such graphs, and $G$ does not contain any of these five as a forbidded subgraph, it follows that $G$ does not contain any of the forbidden subgraphs characterizing \gatex graph, so $G$ is a \gatex graph. Finally, since all induced $P_4$s of $G$ belong to the same connected component $G'$ of $G$, it follows that all connected components of $G$ distinct from $G'$ are cographs.
\end{proof}

\section{Conclusion}\label{sec-out}

In this contribution, we showed that the class of undirected graphs that can be explained by a labelled arboreal network $(N,t)$ is precisely the class of distance-hereditary graphs (Theorem~\ref{th:dh}). As a companion to this result, we provided an algorithm that constructs, for a given distance-hereditary graph $G$, a labelled arboreal network $(N,t)$ explaining $G$ (Algorithm~\ref{alg:dh}).

In addition, we investigated the relationships between the graph explained by a given labelled arboreal network $(N,t)$, and the Ptolemaic graph associated to $N$ (Theorem~\ref{th:ptext}). We also characterized the subclass of distance-hereditary graphs that can be explained by a labelled arboreal network with exactly two roots (Theorem~\ref{th:2rgt}).

These results led to many interesting questions. First of all, many optimization problems can be tackled, for cographs, by algorithms working with the structure of labelled rooted trees, rather than with the structures of cographs themselves. Labelled arboreal networks offer a new playground to generalize these algorithms, and provide new solutions to classical graph problems for distance-hereditary graphs.

Although we know which graphs can be explained by a labelled arboreal network with exactly one root (the cographs), and we characterized here the graphs that can be explained by a labelled arboreal network with exactly two roots (Theorem~\ref{th:2rgt}), we do not know, for a given distance-hereditary graph $G$ that does not fall into one of these two cases, the minimum $k \geq 3$ such that $G$ can be explained with a labelled arboreal network with exactly $k$ roots. In particular, Algorithm~\ref{alg:dh} is not guaranteed to produce a labelled arboreal network with a minimum number of root (see \emph{e. g.} Figure~\ref{fig-alg}). It would thus be of interest to identify this minimum value $k$, and to improve Algorithm~\ref{alg:dh} in such a way that the network returned by the algorithm is optimal in that sense.

Finally, it follows from Theorem~\ref{th:ptext} that, if $(N,t)$ is a labelled arboreal network, then $G^*=\mathcal A(N)$ is a \emph{chordal completion} of $G=\mathcal C(N,t)$, that is, a chordal supergraph of $G$ obtained by adding edges between existing vertices of $G$. We say that a chordal completion $G^*=(X,E^*)$ of a graph $G=(X,E)$ is \emph{minimal} is there exists no $e \in E^* \setminus E$ such that $G'=(X,E^*-\{e\})$ is chordal. We say that $G^*$ is \emph{minimum} if there exists no chordal completion $G'=(X,E')$ of $G$ with $|E'|<|E^*|$. Note that a minimum chordal completion is necessarily minimal. These concepts lead to the following, open question: Given a connected distance-hereditary graph $G$, is it true that, if $G^*$ is a minimum chordal completion of $G$, there exists a labelled arboreal network $(N,t)$ explaining $\mathcal C(N,t)=G$ such that $\mathcal A(N)=G^*$? Interestingly, the answer is no if we replace "minimum" with "minimal".

\begin{acknowledgements}
This contribution finds itself at the intersection of two exciting and still ongoing projects, both of which involving fruitful collaborations. For the "arboreal network" part, I would like to thank Katharina Huber and Vincent Moulton. For the "graph representation" and "\gatex" part, and for his useful suggestions that helped improve the manucript, I would like to thank Marc Hellmuth. Finally, I would like to thank Annachiara Korchmaros, for some pieces of useful advices.
\end{acknowledgements}

\section*{Declarations}

The authors declare that they have no conflict of interest.

\section*{Data availability}

No data were used or generated as part of this work.

\bibliographystyle{spmpsci} 
\bibliography{bibli}

\end{document}